\documentclass[review]{elsarticle}

\usepackage{theorem}
\usepackage{amssymb,amsmath}

\usepackage{mathptmx}       
\usepackage{helvet}         
\usepackage{courier}        
\usepackage{type1cm}        
\usepackage{makeidx}         
\usepackage{graphicx}        
\usepackage{multicol}        
\usepackage[bottom]{footmisc}
\usepackage{subfigure}
\usepackage{booktabs}

\usepackage{epsfig} 
\usepackage{epstopdf}
\usepackage{color}
\usepackage{colordvi}

\newcommand{\duality}[4]         { \big \langle {#1},{#2} \big \rangle_{{#3}\times{#4}} }       
       
\newcommand{\norm}[2]         { \| {#1} \|_{#2} }                      
\newcommand{\bfm}[1]             { \mathbf{#1}     }             %
\newcommand{\qq}             { \bfm{q}     }             
\newcommand{\rr}             { \bfm{r}     }             
\newcommand{\Nabla}       { \boldsymbol{\nabla} }   
\newcommand{\SCZO}            { C^0(\Omega) }                    
\newcommand{\SLTO}            { L^2(\Omega) }                       
\newcommand{\SLTOvec}            { [L^2(\Omega)]^2 }                   
\newcommand{\SHOO}            { H^1(\Omega) }                         
\newcommand{\SHOP}            { H^1(\Ph) }                              
\newcommand{\VVK}            { {V(K_m)} }              
\newcommand{\VV}            { {V(\Ph)} }        
\newcommand{\WW}            { {W(\Ph)} }        
\newcommand{\UUU}            { U(\Omega) }           
\newcommand{\UUUh}            { U^h(\Omega) }           
\newcommand{\VVh}            { V^*(\Ph) }          
\newcommand{\SHOK}            { H^1(K_m) }                        
\newcommand{\SHdivO}            { H(\text{div},\Omega) }             
\newcommand{\SHdivK}            { H(\text{div},K_m) }            
\newcommand{\SHMHGN}            { H^{-1/2}(\GN) }          
\newcommand{\SHHdK}            { H^{1/2}(\dKm) }         
\newcommand{\SHmHdK}            { H^{-1/2}(\dKm)  }               
\newcommand{\SLOinf}            { L^\infty(\Omega)}          
\newcommand{\SLOinft}            { [L^\infty(\Omega)]^2 }   

\newcommand{\DD}            { {\bfm{D}}   }              
\newcommand{\bb}             { {\bfm{b}}   }             
\newcommand{\xx}             { \bfm{x}     }             
\newcommand{\vn}             { \bfm{n}     }             
\newcommand{\err}             { \bfm{\mathcal{E} }    }       
\newcommand{\res}             { \mathcal{R}_h     }       
\newcommand{\eff}             { \mathcal{I}_{eff}     }       

\newcommand{\GD}          { \Gamma_{D} }            
\newcommand{\GN}          { \Gamma_{N} }            
\newcommand{\Ph}            { \mathcal{P}_h }
\newcommand{\Kep}           { K_m \in \Ph}
\newcommand{\dKm}           { \partial K_m }
\newcommand{\dOm}           { \partial \Omega }

\newcommand{\dx}              { \; {\rm d} \bfm{x}   }                
\newcommand{\dss}             { \, {\rm d} s   }                          
\newcommand{\summa}[2]        { \overset{#2}{\underset{#1}{\sum}} } 
\newcommand{\supp}[1]         { \underset{#1}{\sup} \, }        

\newcommand{\Pe}              {  {\rm Pe}   }                
\newcommand{\wwwh}              { \bfm{w^*} }                      
\newcommand{\vvh}              { v^* }                   
\newcommand{\vvi}              { \tilde{e}^i }                 
\newcommand{\wwwi}              { \mathbf{\tilde{E}^i}}          
\newcommand{\vvj}              { \tilde{e}^j_x }           
\newcommand{\wwwj}              { \mathbf{\tilde{E}^j_x}}         
\newcommand{\vvk}              { \tilde{e}^k_y }         
\newcommand{\wwwk}              { \mathbf{\tilde{E}^k_y}}      
\newcommand{\vvih}              { \tilde{e}^i_h }                      
\newcommand{\wwwih}              { \mathbf{\tilde{E}^i_h}}                     
\newcommand{\vvjh}              { \tilde{e}^{j}_{x_h} }                      
\newcommand{\wwwjh}              { \mathbf{\tilde{E}^{j}_{x_h}}}                     
\newcommand{\vvkh}              { \tilde{e}^k_{y_h} }                      
\newcommand{\wwwkh}              { \mathbf{\tilde{E}^k_{y_h}}}                     
\newcommand{\zzz}              { \bfm{z} }                      
\newcommand{\www}             { \bfm{w} }                      
\newcommand{\zz}              { \bfm{z} }                      
\newcommand{\nn}              { \bfm{n} }                      
\newcommand{\isdef}           { \overset{\text{def}}{=} } 
\newcommand{\ds}              { \displaystyle }   

\theoremstyle{plain}
\theoremheaderfont{\bfseries \upshape}

\newtheorem{lem}{Lemma}[section]
\newtheorem{rem}{Remark}[section]

\begin{document}

\begin{frontmatter}
 \title{Goal-Oriented Error Estimation for the Automatic Variationally Stable
FE Method for Convection-Dominated Diffusion Problems}

    \author[1]{Eirik Valseth\corref{cor1}}
\ead{Eirik@utexas.edu}

\author[2]{Albert Romkes}
\ead{Albert.Romkes@sdsmt.edu}

 \address[1]{Oden Institute for Computational Engineering and Sciences, The University of Texas at Austin, Austin, TX 78712, USA}

 \address[2]{Department of Mechanical Engineering, South Dakota School of Mines \& Technology, 501 E. St. Joseph Street, Rapid City, SD 57701, USA}



 \cortext[cor1]{Corresponding author}



\biboptions{sort&compress}


%
%
%
\begin{abstract}
We present goal-oriented \emph{a posteriori} error estimates for the 
automatic variationally stable finite element (AVS-FE) method~\cite{CaloRomkesValseth2018} for scalar-valued 
convection-diffusion problems. The AVS-FE method is a Petrov-Galerkin method in which the  test space is broken, whereas the trial space consists of classical FE basis functions, e.g., $C^0$ or Raviart-Thomas functions. We employ the concept of optimal test functions of the 
discontinuous Petrov-Galerkin (DPG) method by Demkowicz and Gopalakrishnan
\cite{Demkowicz4, Demkowicz2, Demkowicz3, Demkowicz5, Demkowicz6}, leading to unconditionally stable FE approximations.
Remarkably, by using $C^0$ or Raviart-Thomas trial spaces, the optimal discontinuous test functions can be computed in a completely decoupled element-by-element
fashion. 

To establish the error estimators we present two approaches: $i)$
following the classical approach of Becker and Rannacher~\cite{becker_rannacher_2001},
i.e., the dual solution is sought in the (broken) test space, and $ii)$  
introducing an alternative approach in which we seek $C^0$, or Raviart-Thomas, AVS-FE 
approximations of the dual solution by using the underlying 
strong form of the dual boundary value problem (BVP).
Various numerical verifications for 2D convection-dominated diffusion BVPs show that the estimates of the approximation error by the new alternative method are highly accurate, while the classical approach  leads to error estimates of poor quality.
Lastly, we present an algorithm for $h-$adaptive processes  based on control of the 
numerical approximation error via the new alternative approach. Numerical verifications show that the estimator maintains high accuracy as the error converges to zero.
\end{abstract}

\end{frontmatter}

\section{Introduction}
\label{sec:introduction}
In adaptive mesh refinement algorithms,
\emph{a posteriori} error estimation~\cite{ainsworth2011posteriori,oden1989toward} is 
 needed to provide quantified assessments of the 
numerical approximation error, as well as error indicators to 
guide the adaptive process.
Residual-based goal-oriented error estimates have been developed for multiple applications and FE methods, e.g, see~\cite{Oden1,becker_rannacher_2001,giles2002adjoint,ladeveze2013new,
prudhomme1999goal}.
It involves the solution of a dual problem, 
which for Bubnov-Galerkin and Petrov-Galerkin methods~\cite{Bubnov1913,petrov1940application,hughes2012finite,oden2012introduction,
reddy1993introduction} suffers from the same numerical 
instabilities as the primal problem in the presence of convection.
Thus, the FE meshes for the dual problem have to be refined so as to adequately 
capture any boundary or internal layers and thereby avoid any numerical instabilities. 
It makes the classical methods unsuitable for goal-oriented error estimation 
for this class of problems.

However, goal-oriented error estimates have been successfully applied to conditionally stable  FE methods by several authors for convection-diffusion problems~\cite{kuzmin2010goal,cnossen2006aspects,schwegler2016goal,
formaggia2001anisotropic,schwegler2014adaptive}.
In these works, stabilization schemes such as the streamlined-upwind Petrov-Galerkin
(SUPG)~\cite{Brooks1982} method are used to stabilize both the primal and dual problems. The reported effectivity of the estimates varies based on the stabilization chosen, type of error to be estimated (approximation vs. modeling), and the Quantity of interest (QoI) used. 
 In~\cite{schwegler2014adaptive,schwegler2016goal},
Schwegler \emph{et al.} explicitly investigate the stabilization of the dual problem and its influence on the estimate. Stabilized discontinuous Galerkin (DG) methods have also been applied to goal-oriented error estimates for convection-diffusion problems with 
success, we refer to~\cite{mozolevski2015robust} and references therein. 
While these error estimation efforts have been successful, the conditionally stable 
nature of the methods do require a priori analyses to properly establish the 
parameters needed to achieve stability which can be extremely arduous and have to be done on a problem-by-problem basis.

Unconditionally stable FE methods such as
the least squares FE method (LSFEM)~\cite{bochevLeastSquares} and the DPG method~\cite{Demkowicz4, Demkowicz2, Demkowicz3, Demkowicz5, Demkowicz6} resolve the issue of conditional numerical instability.
However, these methods generally use a built-in
\emph{a posteriori} error estimator based on the error in the energy norm induced by its bilinear form to drive adaptivity (see, e.g.,~\cite{Demkowicz2,niemi2013automatically,Demkowicz6,keith2016dpg}).
The simplicity and quality of this type of estimator make it the 
most commonly employed for residual minimization techniques such as the DPG method  and  LSFEM. 
In~\cite{keith2017goal}, Keith \emph{et al.} introduce the 
concept of goal-oriented adaptive mesh refinement for the DPG method. 
However, their goal is not so much to estimate the errors in FE computations but rather the
introduction of a new duality theory that 
is used as a vessel for new adaptive mesh refinement strategies.
In~\cite{chaudhry2014enhancing}, the least squares functional is modified by adding terms 
incorporating QoIs to enhance the quality of the built-in estimator 
used to drive the adaptive mesh refinement. A similar approach is taken by Cai, Ku  \emph{et al.} in~\cite{cai2011goal, ku2010posteriori,ku2011posteriori} to both estimate errors and drive mesh \textcolor{black}{refinements}.

While the error estimation in the aforementioned references has been successful, 
the stabilization efforts required can be arduous by demanding in-depth a priori 
error analyses on a problem-by-problem basis. In addition, to our best knowledge there are no
published results for the DPG method that state the effectivity of goal-oriented
\emph{a posteriori} error estimates for convection-diffusion problems. 
Goal-oriented error estimation in the  \mbox{LSFEM} for convection-dominated 
diffusion problems is less attractive due to the highly diffusive nature of its FE approximations for coarse 
meshes which can lead to estimates with poor accuracy. Our goal is therefore to introduce 
a new framework for the goal-oriented \emph{a posteriori} error estimation for the 
automatically stable AVS-FE method that delivers highly accurate predictions 
of the error in user defined QoIs.

The AVS-FE method
introduced by Calo, Romkes and Valseth~\cite{CaloRomkesValseth2018} provides a 
functional setting to analyze singularly perturbed problems, such as convection-dominated 
diffusion. The AVS-FE method 
is a hybrid between the DPG method~\cite{Demkowicz4, Demkowicz2, Demkowicz3, Demkowicz5, Demkowicz6} and the classical mixed FE methods in the sense that the trial space consists of globally continuous functions, while the  
test space consists of piecewise discontinuous functions. Attractive features of the AVS-FE method are its unconditional numerical 
stability property (regardless of the underlying differential operator),  its highly accurate flux approximations, and the ability to compute optimal test functions element-by-element.

In the following, we limit our focus to  stationary scalar-valued convection-diffusion problems.
 In Section~\ref{sec:avs-fe}, we introduce the model problem of 2D scalar-valued convection-dominated diffusion, used notations, as well as a review of the AVS-FE methodology.
Goal-oriented \emph{a posteriori} error estimates are introduced in Sections~\ref{sec:class_estimates} and~\ref{sec:new_GOAM}. 
In Section~\ref{sec:class_estimates}, we introduce the goal-oriented error estimates for the AVS-FE method following the classical approach of Becker and Rannacher~\cite{becker_rannacher_2001}, i.e., the dual solution is sought in the (broken) test space, and present a numerical verification for the 
Laplace BVP. 
In Section~\ref{sec:new_GOAM}, we present a new alternative approach to goal-oriented 
error estimation in which we seek $C^0$ or Raviart-Thomas AVS-FE 
approximations of the dual solution by using the underlying dual BVP.
Numerical verifications investigating the 
effectivity and robustness of the new estimator are also presented in Section \ref{sec:new_GOAM}. Goal-oriented adaptive mesh refinements and numerical verifications are presented in Section \ref{sec:adaptivity}. Lastly, conclusions and future work are discussed in Section \ref{sec:conclusions}.

\section{Variationally Stable Analysis for Finite Element Computations}  
\label{sec:avs-fe}
In this section, we introduce our convection-diffusion model problem and briefly present a review of the AVS-FE method.
A more detailed introduction can be found in ~\cite{CaloRomkesValseth2018,eirik2019thesis}. 

\subsection{Model Problem and Notation}
\label{sec:model_and_notation}
Let $\Omega \subset \mathbb{R}^2$ be an open bounded domain with Lipschitz boundary $\partial \Omega$ and outward unit normal vector $\vn$. The boundary $\partial \Omega$ consists of  open subsections $\GD,\GN\subset\partial\Omega$, such that $\GD\cap\GN=\emptyset$ and $\partial \Omega = \overline{\GD\cup\GN}$. For our model problem, we  consider the following linear convection-diffusion PDE in $\Omega$ with homogeneous Dirichlet conditions on $\GD$ and (possibly) non-homogeneous Neumann conditions on $\GN$:
\begin{equation} \label{eq:conv_diff_reac_BVP}
\boxed{
\begin{array}{l}
\text{Find }  u  \text{ such that:}    
\\[0.05in] 
\qquad 
\begin{array}{rcl}
  -\Nabla \cdot(\DD \Nabla u) \, + \,
 \bb \cdot \Nabla u  & = & f, \quad \text{ in } \, \Omega, 
 \\[0.05in]
 \qquad u &  = & 0, \quad \text{ on } \, \GD , 
 \\
 \qquad \DD\Nabla u \cdot \vn & = & g, \quad \text{ on } \,  \GN,
 \end{array}
 \end{array}
}
\end{equation}
where $\DD$ denotes the second order diffusion tensor, with symmetric and elliptic coefficients $D_{ij}\in\SLOinf$; $\bb\in\SLOinft$ the convection coefficient;
$f\in\SLTO$ the source function; and $g\in\SHMHGN$ the Neumann data. 

\subsection{The AVS-FE Weak Formulation}
\label{sec:weak_AVS_FORM}

For the sake of brevity, we only mention the few key points here of
the derivation of a weak formulation for the AVS-FE method.
We refer to~\cite{CaloRomkesValseth2018,eirik2019thesis} for a more detailed treatment.
We start by introducing a regular partition 
$\Ph$ of $\Omega$ into elements $K_m$, such that:
\begin{equation}
\notag
\label{eq:domain}
  \Omega = \text{int} ( \bigcup_{\Kep} \overline{K_m} ), \quad K_m \cap K_n, \quad m \ne n.
\end{equation}
The partition $\Ph$ is such that any discontinuities in $D_{ij}$ or $\bb$ are 
restricted to the boundaries of each element $\dKm$. We introduce an auxiliary  flux 
variable  $\qq =\{q_x,q_y\}^T= \DD\Nabla u$, and 
recast~\eqref{eq:conv_diff_reac_BVP} as a system of first-order PDEs:
\begin{equation} \label{eq:conv_diff_BVP_first_order}
\boxed{
\begin{array}{l}
\text{Find }  (u,\qq) \in \SHOO\times\SHdivO \text{ such that:}    
\\[0.05in] 
\qquad 
\begin{array}{rcl}
 \DD\Nabla u - \qq  & =  & 0, \quad \text{ in } \, \Omega, 
  \\
  -\Nabla \cdot \qq \, + \,
 \bb \cdot \Nabla u  & = & f, \quad \text{ in } \, \Omega, 
 \\[0.025in]
 \qquad u &  = & 0, \quad \text{ on } \, \GD , 
 \\
 \qquad \qq \cdot \vn & = & g, \quad \text{ on } \,  \GN. 
 \end{array}
 \end{array}
}
\end{equation}
By weakly enforcing the system of PDEs~\eqref{eq:conv_diff_BVP_first_order} \emph{locally}
 on each element $K_m \in \Ph$, applying Green's identity to the term including
the divergence of $\qq$, applying Dirichlet and Neumann conditions on $\partial K_m \cap \Gamma_D$ and $\partial K_m \cap \Gamma_N$, respectively, and subsequently summing all the local contributions we arrive at the following equivalent global variational formulation:  
\begin{equation} \label{eq:weak_form}
\boxed{
\begin{array}{ll}
\text{Find } (u,\qq) \in U(\Omega) & \hspace{-0.05in} \text{ such that:}
\\[0.05in]
 &  \quad B((u,\qq);(v,\www)) = F(v), \quad \forall (v,\www)\in \VV. 
 \end{array}}
\end{equation}
Here, the bilinear form, $B:\UUU\times\VV\longrightarrow \mathbb{R}$,
and linear functional, $F:\VV\longrightarrow \mathbb{R}$, are defined as follows:
\begin{equation} \label{eq:B_and_F}
\begin{array}{c}
B((u,\qq);(v,\www)) \isdef
\ds \summa{\Kep}{}\biggl\{ \int_{K_m}\biggl[  \, \left(  \DD\Nabla u - \qq \right) \cdot \www_m \, 
+  \, \qq \cdot \Nabla v_m \, + \,
 (\bb \cdot \Nabla u) \, v_m  \biggr] \dx\biggr.
 \\[0.1in]
 \hfill \ds 
 - \oint_{\dKm\setminus \overline{\GD\cup\GN}}  \gamma^m_\nn(\qq) \, \gamma^m_0(v_m) \, \dss \biggr\},
 \\[0.15in]
 F(v) \isdef   \ds \summa{\Kep}{}\left\{\int_{K_m} f\,v_m \dx + \oint_{\dKm\cap\GN} g\, \gamma^m_0(v_m) \dss\right\}, 
 \end{array}
\end{equation}
where the trial and test function spaces, $\UUU$ and $\VV$, are:
\begin{equation}
\label{eq:function_spaces}
\begin{array}{c}
\UUU \isdef \biggl\{ (u,\qq)\in \SHOO\times\SHdivO: \; \gamma_0^m(u)_{|\GD} =0\biggr\},
\\[0.15in]
\VV \isdef \biggl\{ (v,\www)\in \SHOP\times\SLTOvec: \,  \gamma_0^m(v_m)_{|\dKm\cap\GD} =0, \; \forall\Kep\biggr\},
\end{array}
\end{equation}
\textcolor{black}{in which the broken $H^1$ space is defined as:
\begin{equation}
\SHOP \isdef \biggl\{ v\in  \SLTO : \;  v_{|K_m} \in \SHOK, \; \forall\Kep\biggr\},
\end{equation}
}
\textcolor{black}{and} norms $\norm{\cdot}{\UUU}:  \UUU \!\! \longrightarrow \!\! [0,\infty)$ and $\norm{\cdot}{\VV}: \VV\! \! \longrightarrow\! \! [0,\infty)$:
\begin{equation}
\label{eq:broken_norms}
\begin{array}{l}
\ds \norm{(u,\qq)}{\UUU} \isdef \sqrt{\int_{\Omega} \biggl[ \Nabla u \cdot \Nabla u + u^2   + (\Nabla \cdot \qq)^2+\qq \cdot \qq\biggr] \dx }.
\\[0.2in]
\ds   \norm{(v,\www)}{\VV} \isdef \sqrt{\summa{\Kep}{}\int_{K_m} \biggl[  h_m^2 \Nabla v_m \cdot \Nabla v_m + v_m^2   + \www_m \cdot \www_m\biggr] \dx },
 \end{array}
\end{equation}
where $h_m = \text{diam}(K_m)$.
The operators $\gamma^m_0: \SHOK: \longrightarrow \SHHdK$ and $\gamma^m_\nn:\SHdivK \longrightarrow \SHmHdK$ denote the local trace and normal trace operators (e.g., see~\cite{Girault1986}).
Note that we employ an engineering notation convention here by using an integral representation of the boundary integrals rather than that of a duality pairing. 
The variational formulation~\eqref{eq:weak_form} is  essentially a DPG
formulation in which only the space $\VV$ is broken.

\begin{lem}
\label{lem:well_posed_cont}
Let $f\in(\SHOP)'$ and $g\in\SHMHGN$. Then, the weak formulation~\eqref{eq:weak_form} 
is well posed and has a unique solution.
\end{lem} 
\emph{Proof}: \textcolor{black}{We provide only an outline of this proof as similar proofs are available in literature. We first note that the kernel of the underlying convection-diffusion 
differential operator is trivial and introduce an equivalent norm on $\UUU$, the energy norm:}
%
%
%
\begin{equation}
\label{eq:energy_norm}
\norm{(u,\qq)}{B} \isdef \supp{(v,\www)\in \VV\setminus \{(0,\mathbf{0})\}} 
     \frac{|B((u,\qq);(v,\www))|}{\norm{(v,\www)}{\VV}}.
\end{equation}
\textcolor{black}{In the philosophy of the DPG method, we identify an optimal test space, spanned by functions
 that are solutions of a Riesz representation problem:
\begin{equation} \label{eq:abstract_riesz_problem}
\begin{array}{rcll}
\ds \left(\, (\hat{p},\hat{\rr});(v,\www)\, \right)_{\VV} &  \! \! =  \! & B((u,\qq);(v,\www)),& \, \forall (v,\www)\in \VV. 
\end{array}
\end{equation}
Then, $B(\cdot,\cdot)$ satisfies the conditions of the Babu{\v{s}}ka Lax-Milgram Theorem~\cite{babuvska197finite} in terms of the energy norm~\eqref{eq:energy_norm} as the action of  the bilinear form is equivalent to an inner product~\eqref{eq:abstract_riesz_problem}. We refer to the important results of~\cite{carstensen2016breaking}, in the analysis of broken Hilbert spaces and variational formulations. In particular,
 it is shown that broken variational formulations based on differential operators with trivial kernels inherit the stability of their unbroken 
counterparts. }
\newline \noindent ~\qed


\begin{rem} \label{rem:alternative_forms}
It is possible to derive other variational statements in which the 
trial space is continuous and the test space is discontinuous. These will be considered in a forthcoming paper.
\end{rem}

\subsection{AVS-FE Discretization}
\label{sec:discretization}
The AVS-FE method seeks numerical approximations $(u^h,\qq^h)$ of $(u,\qq)$ of the
variational formulation~\eqref{eq:weak_form} by using classical FE bases for 
the trial functions $(u^h,\qq^h)$, i.e., we represent the approximations as linear 
combinations of the trial basis functions \mbox{$(e^i(\xx),(E_x^j(\xx), E^k_y(\xx)))\in\UUUh$}  and their corresponding degrees of freedom:
\begin{equation} \label{eq:FE_sol}
u^h(\xx) = \summa{i=1}{N} u^h_i \, e^i(\xx), 
\quad q^h_x(\xx) = \summa{j=1}{N} q^{h,j}_x \, E_x^j(\xx),
\quad q^h_y(\xx) = \summa{k=1}{N} q^{h,k}_y \, E_y^k(\xx).
\end{equation}
Since the solution space $\UUU$ concerns $\SHOO$ and $\SHdivO$ spaces, the FE 
discretizations can employ classical $\SCZO$ or Raviart-Thomas functions.

The test space $\VV$, however, is discontinuous, 
allowing us to construct piecewise discontinuous \emph{optimal} test functions that yield
unconditionally stable discretizations. These functions are constructed by employing 
the DPG philosophy ~\cite{Demkowicz4, Demkowicz2, Demkowicz3, Demkowicz5, Demkowicz6} in which optimal test functions are defined by  
\emph{global} weak problems. Thus, for the trial functions $e^i(\xx),E_x^j(\xx), \text{and } E^k_y(\xx)$, the corresponding 
global optimal test functions $(\vvi,\wwwi)$, $(\vvj,\wwwj)$, and $(\vvk,\wwwk)$ are the solutions of the following Riesz representation problems~\cite{eirik2019thesis,CaloRomkesValseth2018}, respectively:
\begin{equation}
\label{eq:test_problems}
\begin{array}{rcll}
\ds \left(\, (r,\zzz);(\vvi,\wwwi) \, \right)_\VV &  \! \! =  \! & B(\,(e^i,\mathbf{0});(r,\zzz) \, ),& \, \forall (r,\zz)\in\VV, \quad i=1,\dots, N,
\\[0.15in]
\ds \left(\, (r,\zzz); (\vvj,\wwwj) \, \right)_\VV & \!  \!  =  \! & B(\, (0,(E_x^j,0));(r,\zzz)\, ), & \, \forall (r,\zz)\in\VV, \quad j=1,\dots, N,
\\[0.15in]
\ds \left(\, (r,\zzz); (\vvk,\wwwk) \, \right)_\VV & \!  \!  =  \! & B(\, (0,(0,E_y^k));(r,\zzz) \,), & \, \forall (r,\zz)\in\VV, \quad k=1,\dots, N,
\end{array}
\end{equation}
where $\left(\, (\cdot,\cdot); (\cdot,\cdot) \, \right)_\VV $ denotes the broken 
inner product on $\VV$, defined by:
\begin{equation}
\label{eq:inner_prod_V}
\left(\, (r,\zzz); (v,\www) \, \right)_\VV 
\isdef  \summa{\Kep}{}\int_{K_m} \biggl[  h_m^2 \Nabla r_m \cdot \Nabla v_m + r_m\, v_m   + \zzz_m \cdot \www_m\biggr] \dx.
\end{equation} 
\begin{rem}
Remarkably, the broken nature of the test space $\VV$ allows us to compute numerical approximations of the local restrictions of the optimal test functions in a completely decoupled fashion (see~\cite{eirik2019thesis,CaloRomkesValseth2018} for details).
Thus, we solve local restrictions of~\eqref{eq:test_problems}, e.g.,
\begin{equation}
\label{eq:local_test_problems}
\begin{array}{rcll}
\ds \left(\, (r,\zzz);(\vvih,\wwwih) \, \right)_\VVK & = & B_{|K_m}(\,(e^i,\mathbf{0});(r,\zzz) \, ),& \quad \forall (r,\zz)\in\VVK, 
\\[0.1in]
\ds \left(\, (r,\zzz); (\vvjh,\wwwjh) \, \right)_\VVK & = & B_{|K_m}(\, (0,(E_x^j,0));(r,\zzz)\, ), & \quad \forall (r,\zz)\in\VVK, 
\\[0.1in]
\ds \left(\, (r,\zzz); (\vvkh,\wwwkh) \, \right)_\VVK & = & B_{|K_m}(\, (0,(0,E_y^k));(r,\zzz) \,), & \quad \forall (r,\zz)\in\VVK, 
\end{array}
\end{equation}
where $B_{|K_m}(\cdot;\cdot)$ denotes the restriction of  $B(\cdot;\cdot)$ to the
element $K_m$.
Hence, while the optimal test functions are defined by global 
weak statements, their numerical computation can be performed element-by-element (\textcolor{black}{see, e.g.,~\cite{keith2017discrete} for a detailed discussion on this element-wise assembly process}).
\end{rem}

\begin{rem}
Numerical verifications reveal that the local test functions can be computed by using a degree of approximation that is identical to the degree of approximation of their corresponding trial functions. 
\end{rem}

\begin{rem}
The 
choice of $C^0$ or Raviart-Thomas trial functions has the consequence that the optimal 
test functions have the same support as the trial functions~\cite{eirik2019thesis,CaloRomkesValseth2018}.  
\end{rem}

Finally, we introduce the FE discretization of~\eqref{eq:weak_form} governing the AVS-FE approximation $(u^h,\qq^h)\in\UUU$ of $(u,\qq)$ :
\begin{equation} \label{eq:discrete_form}
\boxed{
\begin{array}{ll}
\text{Find } &  (u^h,\qq^h) \in \UUUh \; \text{ such that:}
\\[0.1in]
 &   B((u^h,\qq^h);(\vvh,\wwwh)) = F(\vvh), \quad \forall (\vvh,\wwwh)\in \VVh, 
 \end{array}}
\end{equation}
where the finite dimensional subspace of test functions $\VVh\subset\VV$ is spanned by the numerical approximations of the test functions  $\{(\vvih,\wwwih)\}_{i=1}^N$, $\{(\vvjh,\wwwjh)\}_{j=1}^N$, and $\{(\vvkh,\wwwkh)\}_{k=1}^N$, as computed from the  test function  problems~\eqref{eq:test_problems} and~\eqref{eq:local_test_problems}.

Since we use the DPG methodology here to construct the optimal test space $\VVh$, the
discrete problem~\eqref{eq:discrete_form} \textcolor{black}{satisfies the conditions of the Babu{\v{s}}ka Lax-Milgram Theorem with continuity and inf-sup constants   
of the continuous problem~\eqref{eq:weak_form} scaled by the continuity constant of a Fortin type operator~\cite{nagaraj2017construction} }. It is therefore unconditionally stable for any choice of mesh parameters $h_m$ and $p_m$.
The corresponding global stiffness matrices are symmetric
and positive definite. 

In the following sections we derive error estimates in terms of user defined QoIs of the 
solution.
The QoIs are represented in terms of continuous linear functionals $Q_i : \UUU \rightarrow \mathbb{R}, i = 1,2,\cdots, N_Q$, for example:
\begin{equation}
\label{eq:QoI}
Q_i(u, \qq) = \frac{1}{|\omega|} \int_{\omega} \, u  \dx, 
\end{equation}
Thus, the goal  is to estimate the error $Q(u-u^h, \qq -\qq^h )$.
We introduce residual based \emph{a posteriori} estimates by taking two
distinctive approaches to the solution of the dual problem. The first follows the 
approach introduced by Becker and Rannacher~\cite{becker_rannacher_2001} and therefore 
seeks a dual solution in the broken primal test space $\VV$.
The second approach concerns an alternative approach in which the AVS-FE solution of the underlying 
strong form of the dual problem is sought 
in a $\SHOO\times\SHdivO$ subspace of the primal test space.

\section{Goal-Oriented Error Estimation - Classical Approach}
\label{sec:class_estimates}
Following Becker and Rannacher~\cite{Oden1,becker_rannacher_2001,
prudhomme1999goal} we state the following classical lemma of goal-oriented 
error estimation:
\begin{lem}
\label{lem:error_in_qoi}
Let $(u,\qq)$ be the exact solution of the first-order system~\eqref{eq:conv_diff_BVP_first_order}, $(u^h,\qq^h)$ the AVS-FE approximation of
$(u,\qq)$ per~\eqref{eq:discrete_form}, and $(p_i, \rr_i) \in \VV $ a dual solution for each QoI, governed by:
\begin{equation} \label{eq:weak_dual}
\boxed{
\begin{array}{ll}
\text{Find }( p_i, \rr_i) \in \VV & \hspace{-0.05in} \text{ such that:}
\\[0.05in]
 &  \quad B((v,\www);(p_i, \rr_i)) = Q_i(v, \www), \quad \forall (v,\www)\in \UUU. 
 \end{array}}
\end{equation}
Then, the error in the QoI $ \err_i(u^h,\qq^h) =  Q_i(u-u^h, \qq-\qq^h)$, is governed by the identity:
\begin{equation} \label{eq:error_identity}
\boxed{
\err_i(u^h,\qq^h) \isdef \res ((u^h,\qq^h);(p_i, \rr_i)),
}
\end{equation}
where $\res ((\cdot,\cdot);(\cdot,\cdot))$ is the residual functional:
\begin{equation} \label{eq:residual_functional}
\res ((u,\qq);(v, \www)) = F(v) - B((u,\qq);(v, \www)).
\end{equation}
\end{lem}
\begin{rem} \label{rem:well_posed_dual}
Analogous to the well-posedness of the  primal problem~\eqref{eq:weak_form} (see Lemma~\ref{lem:well_posed_cont}), the dual problems~\eqref{eq:weak_dual} are well posed (since the kernel of the 
adjoint operator of $B(\cdot,\cdot)$ is also trivial).
\end{rem}
To compute estimates of the error $\err_i(u^h,\qq^h)$ through~\eqref{eq:error_identity}, we compute approximations of
the dual solutions $(p^h_i, \rr^h_i)$
 by following the classical approach of~\cite{Oden1,becker_rannacher_2001,prudhomme2003computable}. Thus, for a given QoI,  the approximate dual solution $(p^h, \rr^h)$ is
governed by:
\begin{equation} \label{eq:dual_discrete_form}
\boxed{
\begin{array}{ll}
\text{Find } &  (p^h, \rr^h) \in \VVh \; \text{ such that:}
\\[0.1in]
 &   B((v^h,\www^h);(p^h, \rr^h)) = Q(v^h,\www^h), \quad \forall (v^h,\www^h) \in \UUUh. 
 \end{array}}
\end{equation}
We seek $(p^h, \rr^h)$ in the discrete, broken space $\VVh$ spanned by $(\vvih,\wwwih)$ of the local Riesz representation problems~\eqref{eq:local_test_problems}.
Hence, we use the same element partition $\Ph$ of $\Omega$ as we used for the primal problem to compute $(p^h, \rr^h)$. However,
due to the Galerkin orthogonality condition of the numerical approximation error, the approximate dual solution has to be sought 
by using polynomial approximations that are of higher order than the trial functions used to solve the discrete primal problem~\eqref{eq:discrete_form}. We choose $p+1$. 
Hence, the 
unconditional numerical stability of the AVS-FE methodology will allow
the computation of approximate dual solutions $(p^h_i, \rr^h_i)$ for any choice of mesh parameters $h_m$ and $p_m$.
The estimated error $\eta_{est}$ in the quantity of interest is then computed by:
\begin{equation} \label{eq:error_estimate}
\eta_{est} \approx \err^h_i(u^h,\qq^h) = \res ((u^h,\qq^h);(p^h_i, \rr^h_i))
\end{equation}
This classical approach has been shown to be very successful in a wide range of
applications,
especially those in which the differential operator is self-adjoint (e.g., see ~\cite{Oden1,prudhomme2003computable}). 

As a numerical verification of this estimator, we consider the Laplace problem
on the unit square $\Omega = (0,1)\times(0,1)\subset \mathbb{R}^2$ with homogeneous Dirichlet boundary conditions:
\begin{equation}\label{eq:poisson_pde}  
\begin{array}{rl}
\ds  - \Delta \, u  = f, & \quad \text{ in } \Omega, 
 \\[0.1in]
  u = 0, & \quad \text{ on } \partial \Omega.  
 \end{array} 
\end{equation}
The source function $f$, is chosen such that the exact solution is given by:
\begin{equation*} 
u(x,y) = e^{\left[ 50(x^2-x)(y^2-y)\right]}-1.
\end{equation*}
The QoI is chosen to be the average of the solution $u$ in the region $\omega = (0.5,1)\times(0.5,1)\subset \Omega$: 
\begin{equation}
\label{eq:QoI_example_1}
Q(u, \qq) = \frac{1}{|\omega|} \int_{\omega} \, u  \dx.
\end{equation}

To estimate the error in this QoI~\eqref{eq:QoI_example_1}  we apply the 
AVS-FE discretization to the primal and dual problem with polynomial degrees of 
approximation of 2 and 3, respectively. As in~\cite{CaloRomkesValseth2018}, we use $C^0$ 
continuous bases for both trial variables $u^h$ and $\qq^h$ of the same polynomial degree,
while for $p^h$ and $\rr^h$ we use the optimal bases determined by the Riesz representation problems~\eqref{eq:test_problems}.
\textcolor{black}{The mesh initial mesh partition used consists of a single quadrilateral element. Subsequent meshes are a sequence uniformly refined from the initial single element.}
\begin{table}[h]
\centering
\caption{\label{tab:classical_results_poisson}  Error estimation results for the Laplace problem with QoI~\eqref{eq:QoI_example_1} using the classical approach, i.e., through~\eqref{eq:dual_discrete_form} and~\eqref{eq:error_estimate}.}
\begin{tabular}{@{}lllll@{}}
\toprule
{Primal dofs \hspace{6mm}} & {$Q(u, \qq)-Q(u^h, \qq^h) $ \hspace{6mm}} & {Dual dofs \hspace{6mm}} & {$ \! \eta_{est} \!$ \hspace{15mm}} & {$\eff \!$} \\
\midrule \midrule

27 & -9.2601e+00 & 48 & -1.0601e+01 & 1.145   \\
75 & 2.3191e-02 & 147 & -7.5337e-02 & -0.325   \\
243 & 1.8610e-02 & 507 & 3.4581e-02 & 1.858   \\
867 & -3.7844e-05 & 1875 & -3.8616e-04 & 10.204   \\
3267 & -2.5778e-05 & 7203 & -5.3075e-05 & 2.059   \\

\bottomrule
\end{tabular}
\end{table}
\begin{figure}[h]
\subfigure[ \label{fig:approx_discont_dual_poisson_a}Dual solution in the classical approach. ]{\centering
 \includegraphics[width=0.5\textwidth]{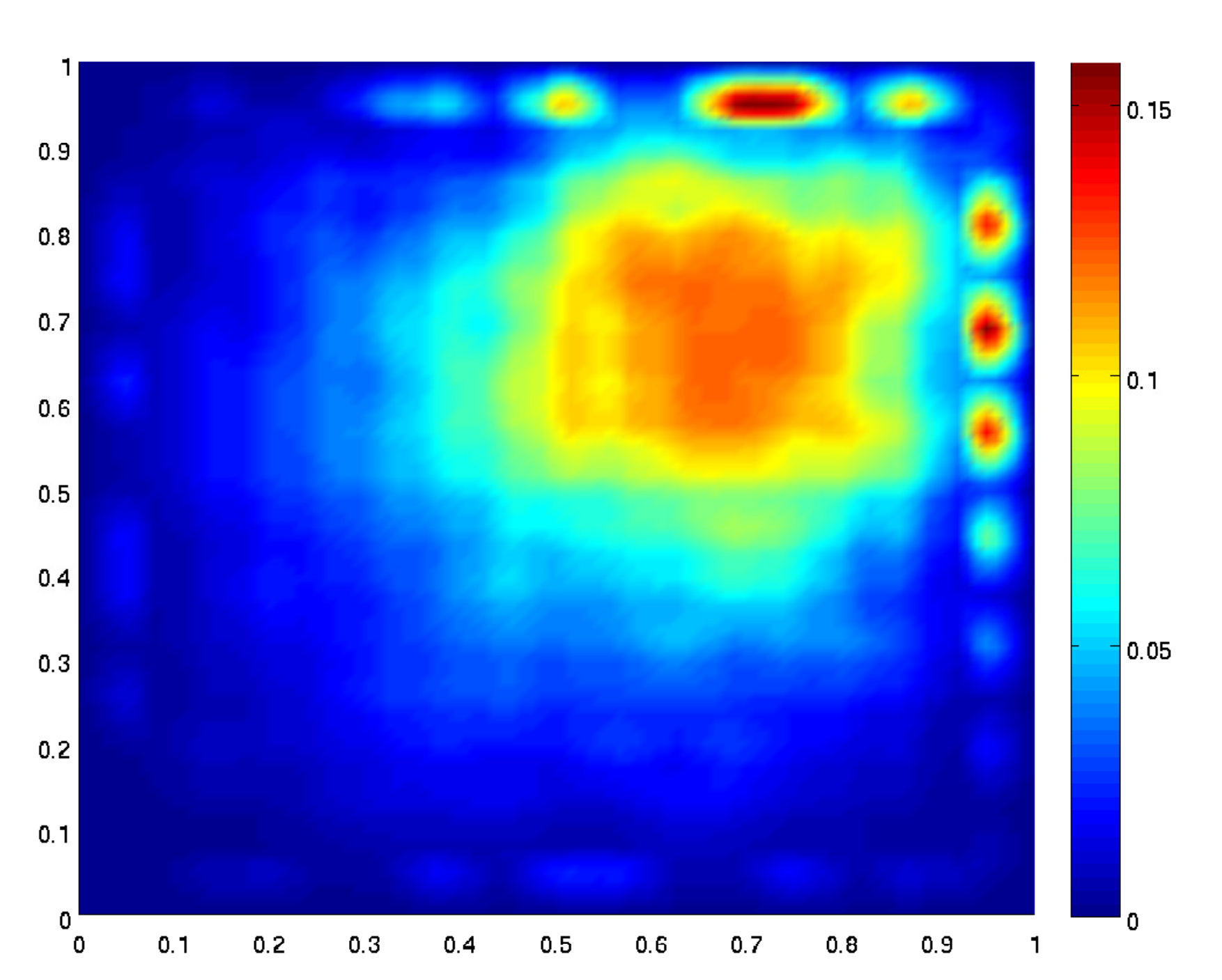}}
\hfill    \subfigure[\label{fig:overkill_dual_poisson} Overkill solution.]{\centering
 \includegraphics[width=0.4875\textwidth]{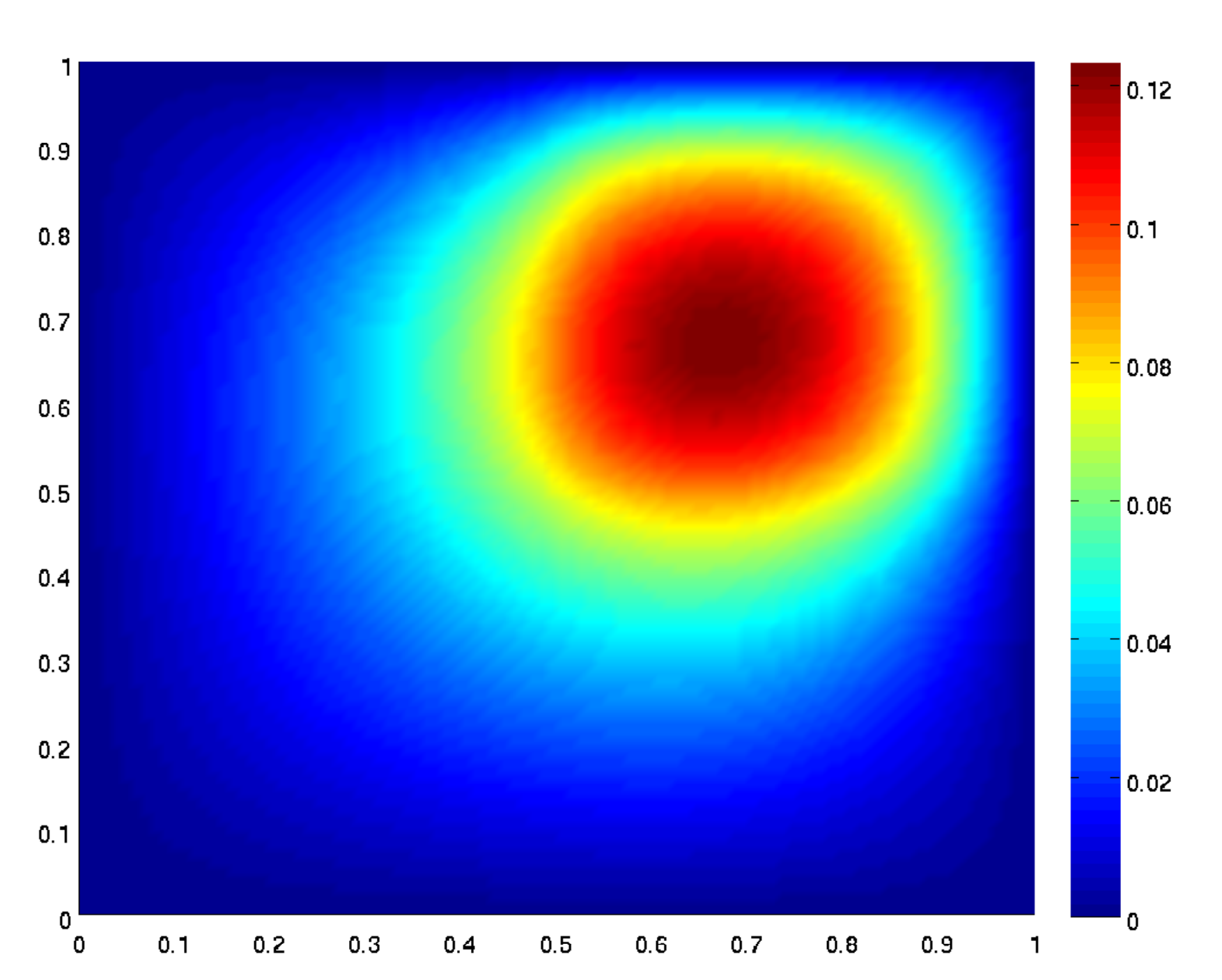}}
  \caption{\label{fig:approx_discont_dual_poisson} Dual solution $p^h$.}
\end{figure}
To assess the quality of the error estimate, we introduce the effectivity index:
\begin{equation} \label{eq:effectivity}
\eff = \frac{\eta_{est}}{Q(u,\qq) - Q(u^h,\qq^h)}.
\end{equation} 
In Table~\ref{tab:classical_results_poisson}, we present the error 
estimates for 
increasingly refined meshes.
It appears that the magnitude of the estimated error decreases monotonically.
However, the results for the effectivity index $\eff$ reveal that the estimates
generally have 
poor accuracy nor exhibit any consistent evolution during the uniform $h-$refinements.
We suspect that since the continuity of the dual solution is enforced weakly in~\eqref{eq:weak_dual}, it prohibits from adequately resolving the discretization of 
the dual solution on the used meshes. \textcolor{black}{We suspect that without drastically changing the formulation using
 jump and/or average operators on the mesh skeleton as one would in, e.g, discontinuous Galerkin methods the discontinuous solutions will not provide adequate resolution of the dual solution}.
Comparison of the approximate solution $p^h$ in Figure~\ref{fig:approx_discont_dual_poisson} to an overkill approximation of $p$ 
in Figure~\ref{fig:overkill_dual_poisson}, reveals internal oscillations in $p^h$  in each element as well as along the global
boundaries. \textcolor{black}{The overkill solution is obtained using a Galerkin FE method for the underlying known dual BVP on a highly refined uniform mesh of quadratic quadrilateral elements.}
While bounded, these contribute to the poor quality of the error estimates. 
\textcolor{black}{These observed oscillations show no consistent behavior during mesh refinements.}
This  behavior persists for numerical verifications in which the degree of approximation for the dual problem is $p_{dual} = p_{primal}+2,p_{primal}+3$ etc.
\textcolor{black}{Another factor that plays a role here can be deduced from DPG$^*$ techniques~\cite{demkowicz2020dpg} since the dual problem~\eqref{eq:dual_discrete_form} can be interpreted as a modified DPG$^*$ form. In particular, it is shown in~\cite{demkowicz2020dpg} that the regularity 
of the domain $\Omega$ is a crucial factor in the accuracy of the dual solution. }

\section{Goal-Oriented Error Estimation - Alternative Approach}
\label{sec:new_GOAM}
Since we suspect that the poor accuracy of the estimator via the classical approach is likely caused by the discontinuous character of the numerical
approximation of the dual solutions $(p^h,\rr^h)$, we propose an alternative approach for computing $(p^h,\rr^h)$. 
Instead of seeking discontinuous discrete approximations of the dual solution by
 using the corresponding dual weak formulation~\eqref{eq:weak_dual} of the primal 
 problem~\eqref{eq:weak_form}, we rather reconsider the underlying strong form
 of each dual problem, i.e., 
\begin{equation} \label{eq:dual_conv_diff_BVP_first_order}
\boxed{
\begin{array}{l}
\text{Find }  (p_i,\rr_i) \in \SHOO\times\SHdivO \text{ such that:}    
\\[0.05in] 
\qquad 
\begin{array}{rcl}
  -\rr_i\, + \Nabla p_i \, & =  & \mathbf{0}, \quad \text{ in } \, \Omega,   \\
  -\Nabla \cdot (\DD \, \rr_i) \, - \,
 \bb \cdot \Nabla p_i   & = & \theta_i, \quad \text{ in } \, \Omega, 
 \\[0.025in]
 \qquad p_i &  = & 0, \quad \text{ on } \, \dOm,
 \end{array}
 \end{array}
}
\end{equation}
where $\theta_i \in U'(\Omega)$ is such that
 $\duality{\theta_i}{(u,\qq)}{U'}{U} = Q_i(u,\qq)$.
\textcolor{black}{Comparison of~\eqref{eq:dual_conv_diff_BVP_first_order} and the first order system for the primal problem~\eqref{eq:conv_diff_BVP_first_order} shows that the diffusion tensor has shifted and the signs in the vector valued PDE has changed. The reason is that~\eqref{eq:dual_conv_diff_BVP_first_order} is the natural form of the distributional first order dual BVP when derived from the AVS-FE weak formulation.}
We subsequently derive a weak statement governing $(p_i,\rr_i)$ by using the same approach 
as the applied to the derivation of the weak statement of the primal problem (see Section~\ref{sec:avs-fe}).
Thus, we seek $p_i \in \SHOO$ , $\rr_i \in \SHdivO$ and employ test spaces for the dual 
problem that are broken. Hence, $p_i$ and $\rr_i$  belong to  the same globally (weakly) continuous function spaces as the primal solution 
$(u,\qq)$. To derive the dual weak statement, we follow the derivation in~\cite{CaloRomkesValseth2018} and enforce the system~\eqref{eq:dual_conv_diff_BVP_first_order} weakly on each element $\Kep$, 
apply Green's Identity, enforce boundary conditions, and arrive at the following weak statement:
\begin{equation} \label{eq:weak_form_dual_new}
\boxed{
\begin{array}{ll}
\text{Find } (p_i,\rr_i) \in U(\Omega) & \hspace{-0.05in} \text{ such that:}
\\[0.05in]
 &  \quad \hat{B}((v,\www);(p_i,\rr_i)) = Q_i(v,\www), \quad \forall (v,\www)\in \WW,
 \end{array}}
\end{equation}
where:
\begin{equation} \label{eq:dual_hat_B}
\begin{array}{l}
\ds \hat{B}((v,\www);(p_i,\rr_i)) = \summa{\Kep}{}\biggl\{ \int_{K_m}\biggl[  \, \left( \Nabla p_i - \rr_i \right) \cdot \www_m \, 
+  \, \DD \, \rr_i \cdot \Nabla v_m \, - \,
 (\bb \cdot \Nabla p_i) \, v_m   \biggr] \dx\biggr.
 \\[0.1in]
 \hspace{1in} \biggl. \ds - \oint_{\dKm\setminus \overline{\dOm}}  \gamma^m_\nn(\DD \, \rr_i) \, \gamma^m_0(v_m) \, \dss \biggr\},
 \end{array}
\end{equation}
and:
\begin{equation}
\label{eq:function_space_W}
\begin{array}{c}
\WW \isdef \biggl\{ (v,\www)\in \SHOP\times\SLTOvec: \,  \gamma_0^m(v_m)_{|\dKm\cap\dOm} =0, \; \forall\Kep\biggr\}.
\end{array}
\end{equation}

It should be noted here that $\hat{B}(\cdot;\cdot)$ and $B(\cdot;\cdot)$ 
differ in the sign in front of the convection vector $\bb$ due to the non-self adjoint
character of the differential operator of~\eqref{eq:conv_diff_reac_BVP}. Now,
to ensure the unconditional stability of the discrete dual problem 
we use optimal discontinuous test functions $(\vvi,\wwwi)$, $(\vvj,\wwwj)$, and $(\vvk,\wwwk)$ for the dual problem that are solutions of the following (Riesz) weak problems:
\begin{equation}
\label{eq:dual_test_problems}
\begin{array}{rcll}
\ds \left(\, (x,\zzz);(\vvi,\wwwi) \, \right)_\VV &  \! \! =  \! & \hat{B}(\,(x,\zzz);(e^i,\mathbf{0}) \, ),& \, \forall (x,\zz)\in\VV, \quad i=1,\dots, N,
\\[0.15in]
\ds \left(\, (x,\zzz); (\vvj,\wwwj) \, \right)_\VV & \!  \!  =  \! & \hat{B}(\, (x,\zzz);(0,(E_x^j,0))\, ), & \, \forall (x,\zz)\in\VV, \quad j=1,\dots, N,
\\[0.15in]
\ds \left(\, (x,\zzz); (\vvk,\wwwk) \, \right)_\VV & \!  \!  =  \! & \hat{B}(\, (x,\zzz);(0,(0,E_y^k)) \,), & \, \forall (x,\zz)\in\VV, \quad k=1,\dots, N.
\end{array}
\end{equation}
\begin{rem}
The vector valued dual solution $\DD\rr_i$ belongs to $\SHdivO$ due to the boundary integral
$\oint_{\dKm\setminus \overline{\dOm}}  \gamma^m_\nn(\DD \, \rr_i) \, \gamma^m_0(v_m) \, \dss$ 
in~\eqref{eq:weak_form_dual_new}. 
For this integral to be Lebesgue integrable, $\DD \, \rr_i$ has to belong to $\SHmHdK$ which implies $\DD \, \rr_i \in\SHdivO $.
\end{rem}
We then establish the error estimator $\hat{\eta}_{est}$ by using the new dual solutions: 
\begin{equation} \label{eq:error_estimate2}
\ds \hat{\eta}_{est} \approx \err^h_i(u^h,\qq^h) = \res ((u^h,\qq^h);(p^h_i, \rr^h_i))
\end{equation}

Having established the new alternative error estimates, we 
propose to employ an error indicator $\epsilon_m$ corresponding to the restriction of the goal-oriented error estimate $\hat{\eta}_{est}$ in mesh adaptive 
refinements, i.e,
\begin{equation} \label{eq:error_indicator}
\epsilon_m = \res|_{K_m} ((u^h,\qq^h);(p^h_i, \rr^h_i)).
\end{equation}

\textcolor{black}{
\begin{rem}
We note that the philosophy we advocate here 
for the dual problem was also proposed  for the consideration of adjoint equations in inverse FE methods  by considering the strong form 
of the adjoint equation  by Bramwell in~\cite{bramwell2013discontinuous}. 
\end{rem}}

\subsection{Numerical Verification - Diffusion Problem}
\label{sec:alternative-app-example-poisson}
We again solve~\eqref{eq:poisson_pde} using
identical meshes and degrees of approximation as in Section~\ref{sec:class_estimates}, i.e., quadratic primal and cubic dual approximations, respectively. In this alternative approach, we 
seek $C^0$ continuous solutions to both the primal and dual problems in which the 
scalar and flux variables are of the same polynomial degree. Again, we assess the quality of the error estimate with the effectivity index:
\begin{equation} \label{eq:effectivity_new}
\hat{\eff} = \frac{\hat{\eta}_{est}}{Q(u,\qq) - Q(u^h,\qq^h)}.
\end{equation} 
\begin{table}[h]
\centering
\caption{\label{tab:new_results_poisson} Error estimation results for the Laplace problem with QoI~\eqref{eq:QoI_example_1} using the new alternative approach, i.e., through~\eqref{eq:weak_form_dual_new} and~\eqref{eq:error_estimate}.}
\begin{tabular}{@{}lllll@{}}
\toprule
{Primal dofs \hspace{6mm}} & {$Q(u, \qq)-Q(u^h, \qq^h) $ \hspace{6mm}} & {Dual dofs \hspace{6mm}} & {$ \! \hat{\eta}_{est} \!$ \hspace{15mm}} & {$\hat{\eff} \!$} \\
\midrule \midrule

\textcolor{black}{27} & \textcolor{black}{-9.2601e+00} & \textcolor{black}{48} & \textcolor{black}{-1.0601e+01} & \textcolor{black}{1.145}   \\
\textcolor{black}{75} & \textcolor{black}{2.3192e-02} & \textcolor{black}{147} & \textcolor{black}{3.4410e-02} & \textcolor{black}{1.484}   \\ 
243 & 1.8610e-02 & 507 & 1.8602e-02 & \textcolor{black}{0.999}   \\
867 & -3.7845e-05 & 1875 & -3.6160e-05 & \textcolor{black}{0.956}   \\
3267 & -2.5778e-05 & 7203 & -2.3479e-05 & \textcolor{black}{0.911}  \\

\bottomrule
\end{tabular}
\end{table}
\begin{figure}[h]
\subfigure[ Surface view. ]{\centering
 \includegraphics[width=0.5\textwidth]{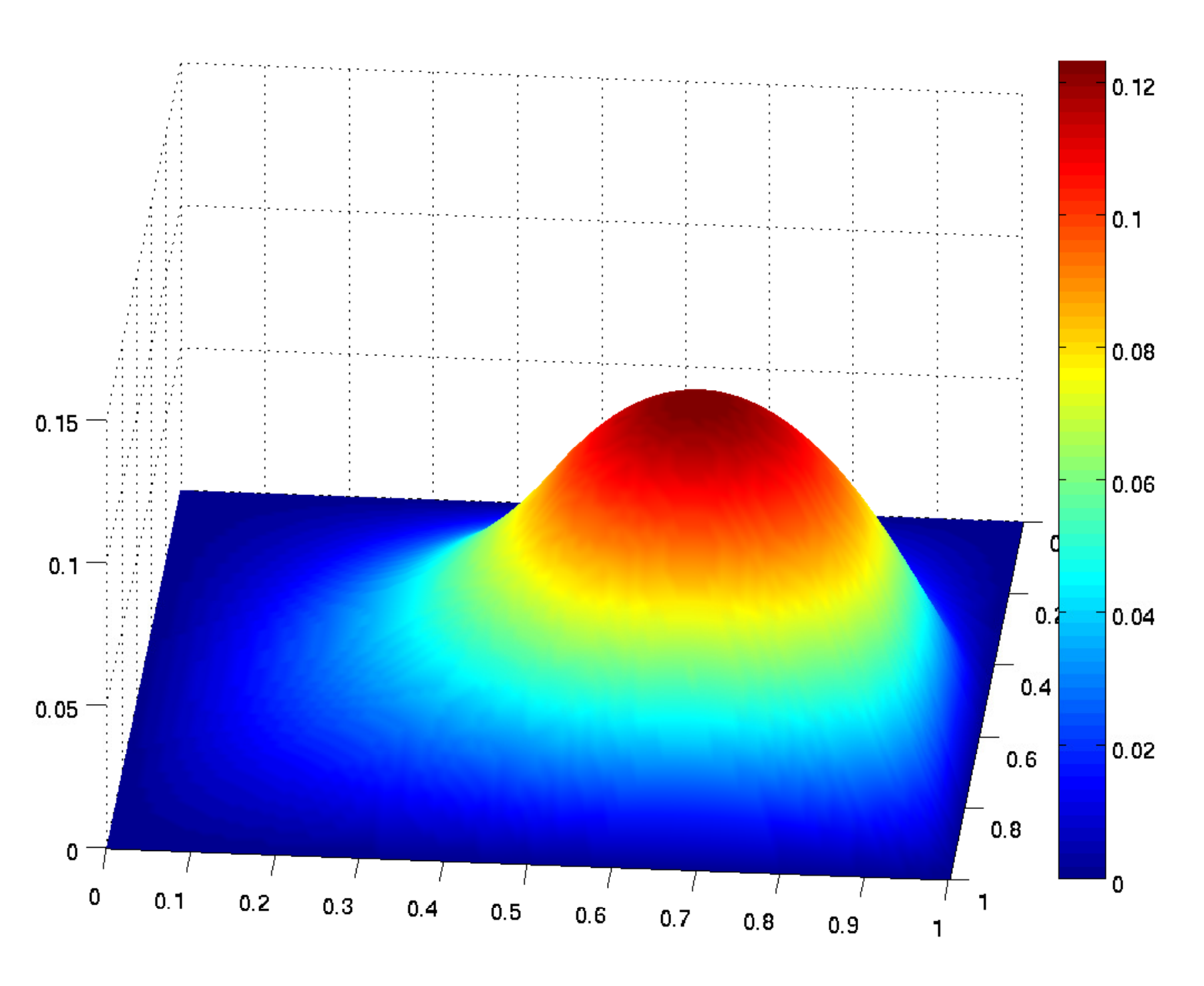}}
\hfill    \subfigure[ Top view.]{\centering
 \includegraphics[width=0.4875\textwidth]{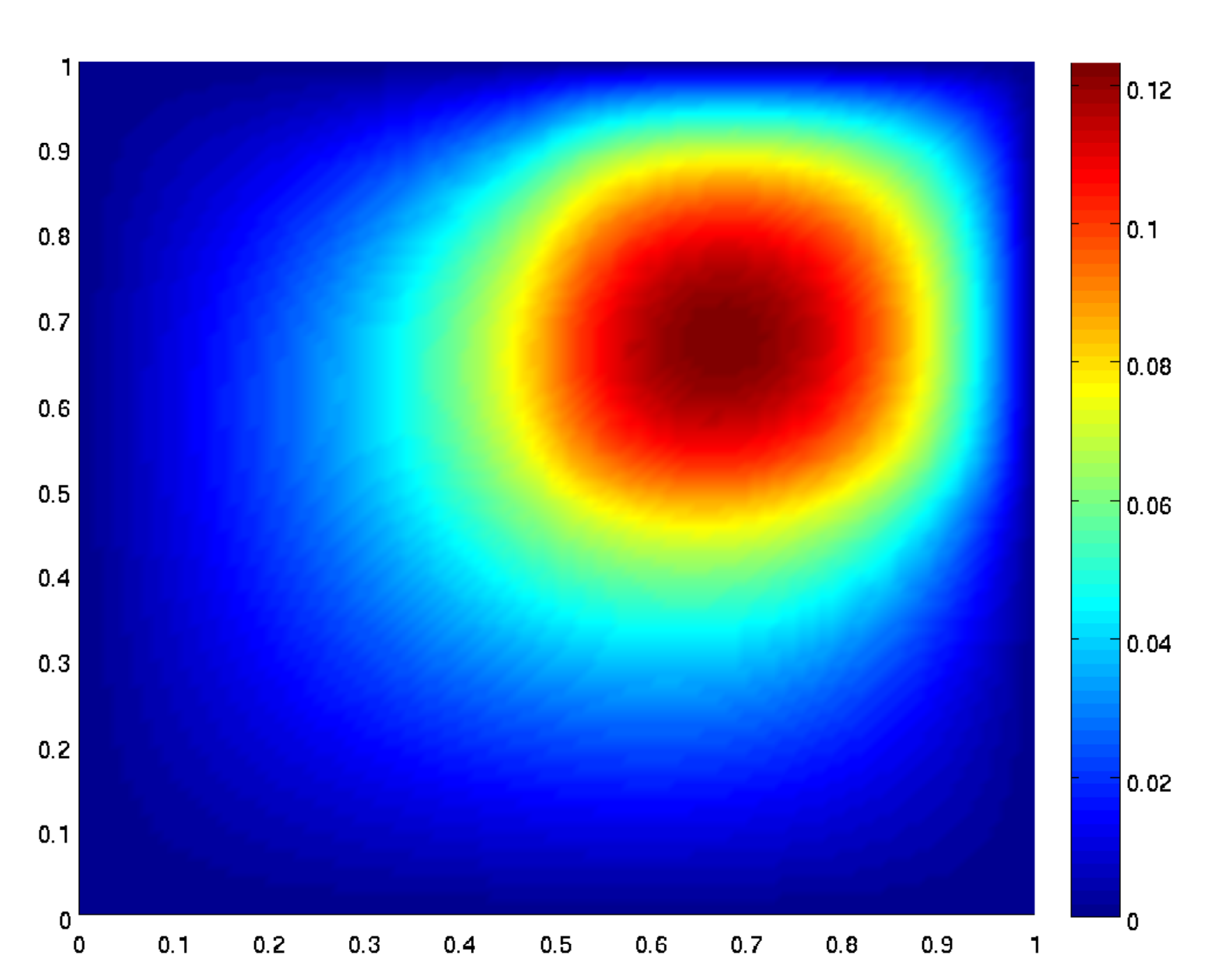}}
  \caption{\label{fig:approx_cont_dual_poisson} $C^0$ dual solution $p^h$ obtained by the new alternative approach through~\eqref{eq:weak_form_dual_new} on a uniform mesh with 3267 dofs.}
\end{figure}
In Table~\ref{tab:new_results_poisson}, we show the results for the error estimates for 
uniform mesh refinements.
As in Section~\ref{sec:class_estimates}, the 
magnitude of the error decreases monotonically. However, the effectivity index now 
shows that the estimates have very good accuracy with values close to unity.
Furthermore, comparison of Figure~\ref{fig:approx_cont_dual_poisson} of the dual 
AVS-FE solution to an overkill solution reveals that there are no oscillations at 
element interiors and global boundaries.

\subsection{Numerical Verifications - Convection-Dominated Diffusion}
\label{sec:conv-diff-ex} 
\begin{figure}[h]
{\centering
 \includegraphics[width=0.8\textwidth]{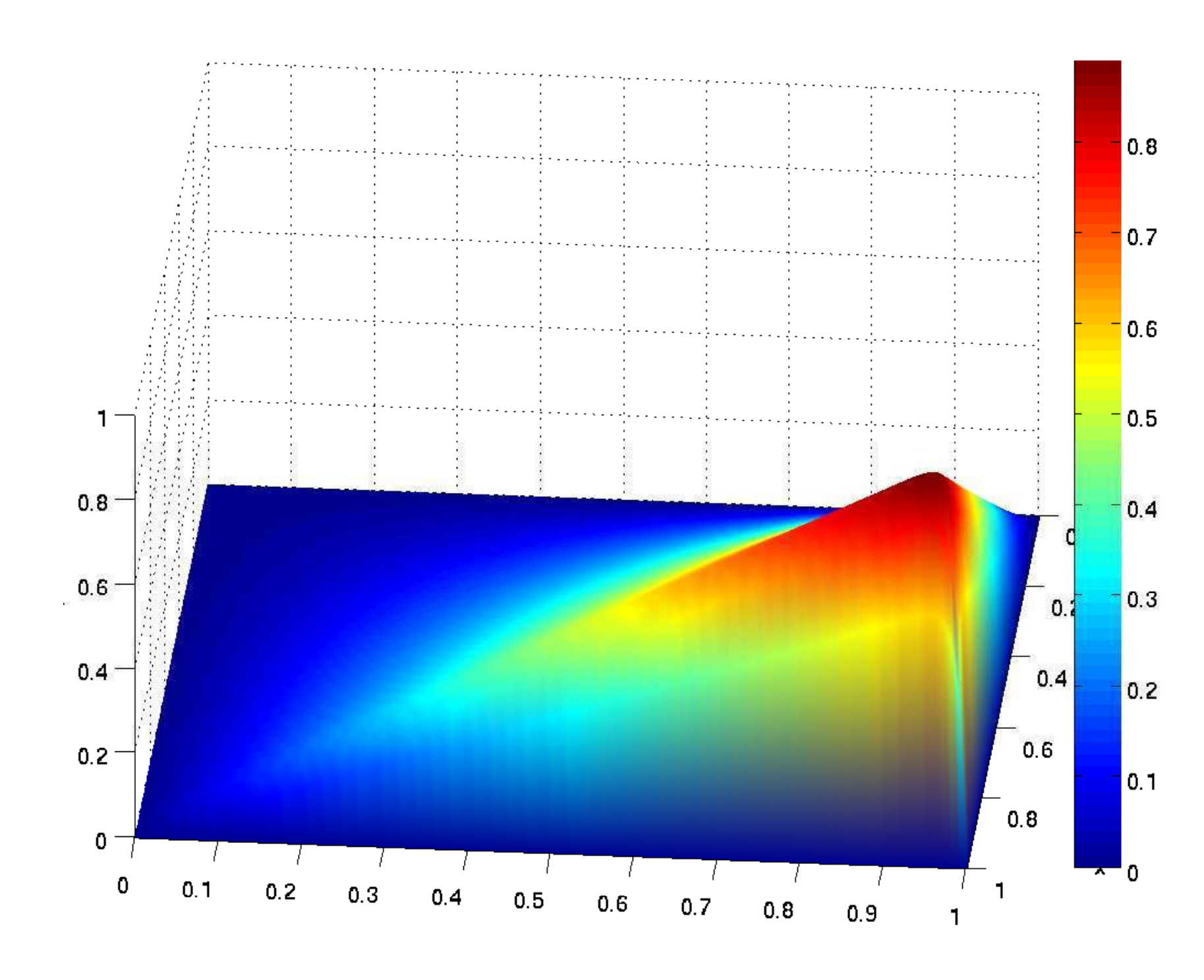}}
 \caption{\label{fig:exact_Pe_100} Exact solution $u$ of the simplified model problem~\eqref{eq:conv_diff_ex} with $\bb=\{1,1\}^T$ and $\Pe=100$. }
\end{figure}

Next, we consider a more challenging case of a convection-dominated diffusion problem. In this subsection,  
we consider a simplified 
form of our model problem~\eqref{eq:conv_diff_reac_BVP} on the unit square with homogeneous Dirichlet boundary conditions:
\begin{equation}\label{eq:conv_diff_ex}  
\begin{array}{rl}
\ds  - \frac{1}{\Pe} \Delta u +
 \bb \cdot  \Nabla u = f, & \quad \text{ in } \Omega = (0,1)\times(0,1), 
 \\[0.1in]
  u = 0, & \quad \text{ on } \partial \Omega,  
 \end{array} 
\end{equation}
where the P\' eclet number $\Pe = 100$ is and $\bb=\{1,1\}^T$ the convection coefficient.
We consider the case of~\eqref{eq:conv_diff_ex} in which 
the above source function $f$ is chosen such that the exact solution is  given by:
\begin{equation} \label{eq:conv_diff_exact}
u(x,y) = \left[x + \frac{e^{\Pe \cdot b_x \cdot x}-1}{1-e^{\Pe \cdot b_x}}\right]\left[y + \frac{e^{\Pe \cdot b_y \cdot y}-1}{1-e^{\Pe \cdot b_y}}\right].
\end{equation}
Thus, the solution 
exhibits boundary layers along $x=1$ and $y=1$ with a width of $\frac{1}{\Pe} = 0.01$, as shown in Figure~\ref{fig:exact_Pe_100}.

\subsubsection{Uniform Meshes}
\label{sec:uniform_mesh}
\begin{figure}[h]
\subfigure[ Surface view. ]{\centering
 \includegraphics[width=0.5\textwidth]{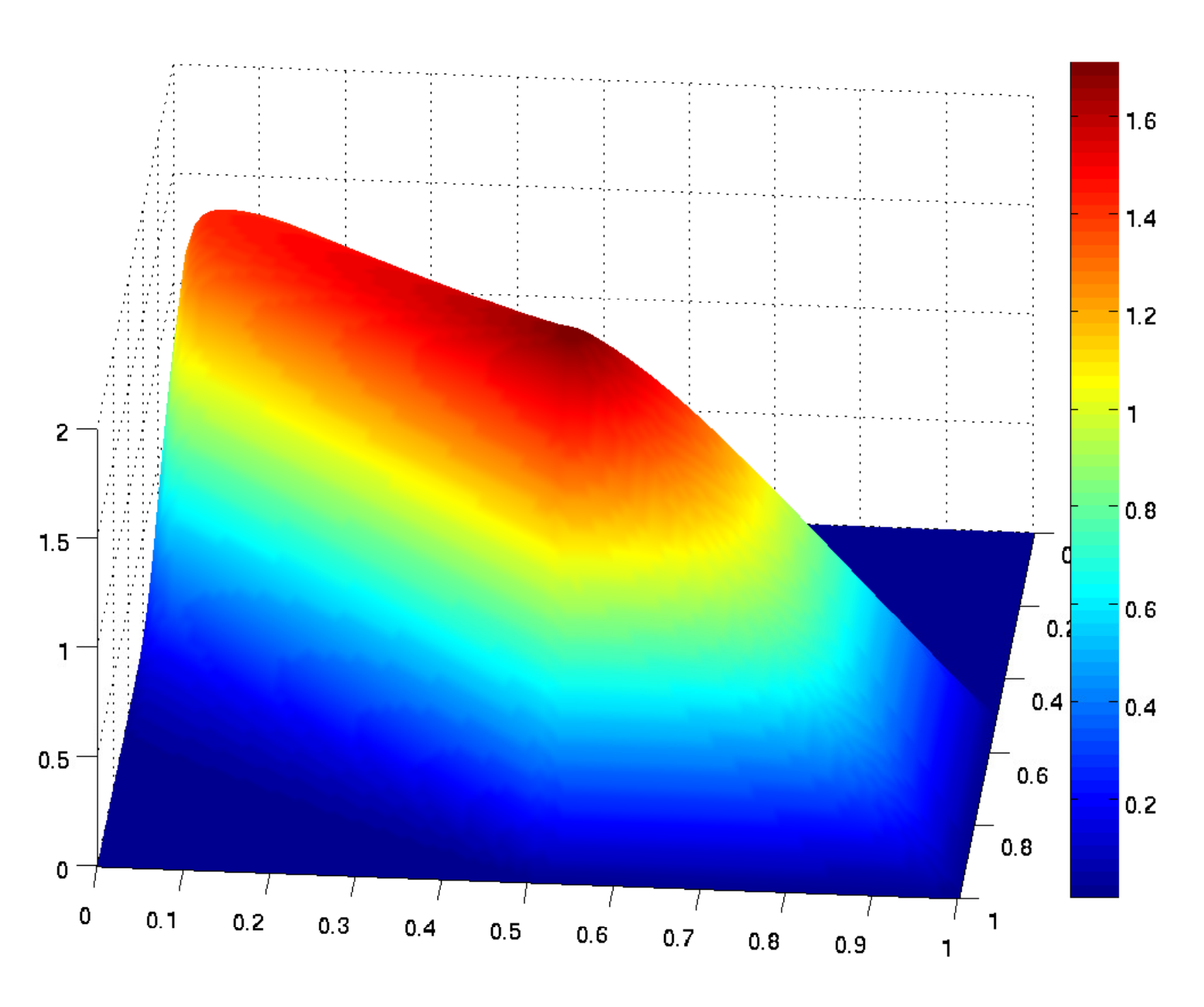}}
\hfill    \subfigure[ Top view.]{\centering
 \includegraphics[width=0.4875\textwidth]{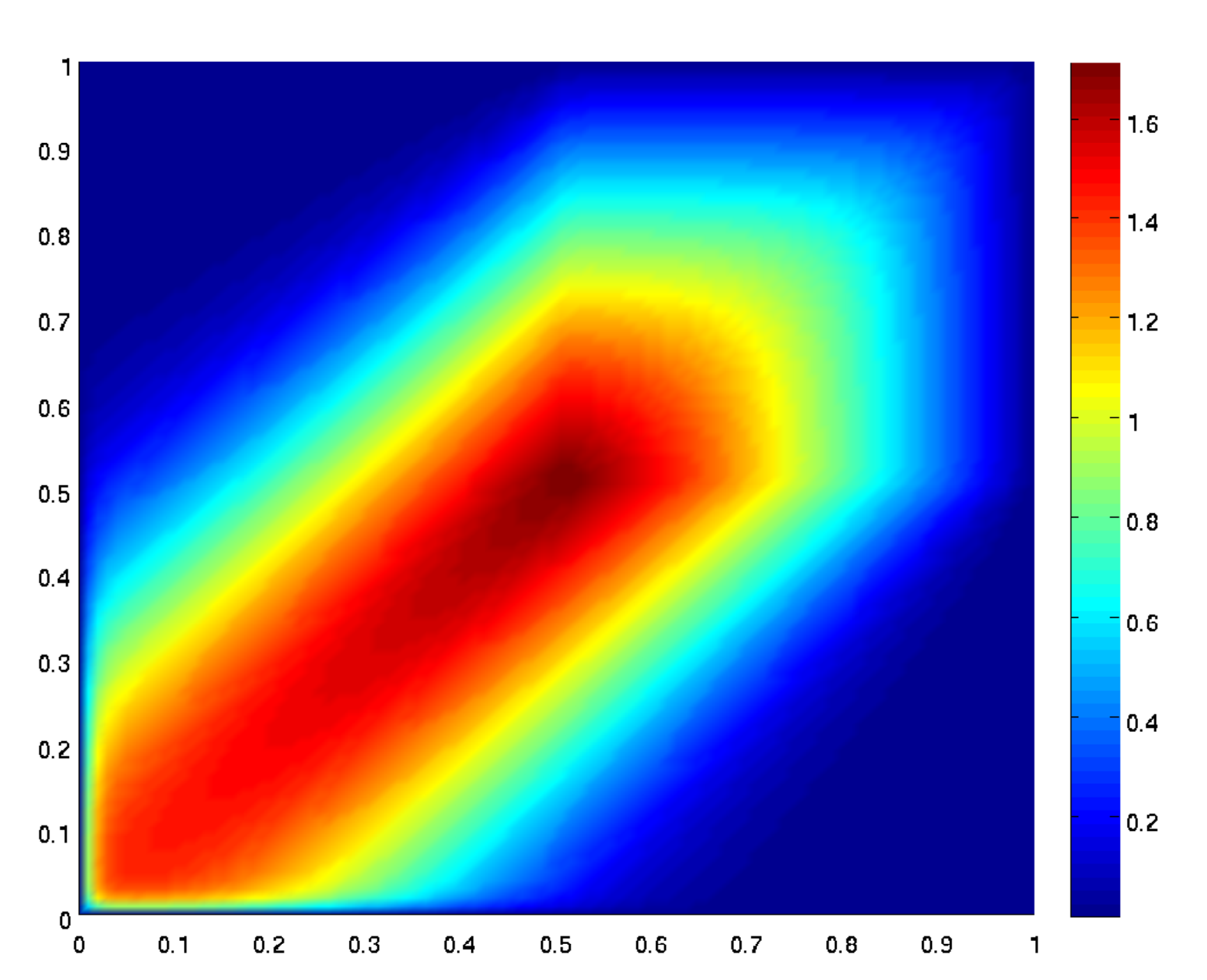}}
  \caption{\label{fig:dual_overkill} Overkill solution $p$ of the dual problem with $\bb=\{1,1\}^T$ and $\Pe=100$.}
\end{figure}

First, we consider uniform meshes consisting of quadrilateral elements.
The first QoI is chosen as in~\eqref{eq:QoI_example_1}, i.e., the average solution $u$ in
the top right quadrant of the unit square.
In Figure~\ref{fig:dual_overkill}, we show the corresponding overkill solution $p$ of
the dual problem.
The dual solution, similar to the primal, exhibits boundary layers. However, as the direction 
of the convection is reversed from the primal problem,
 the layers are at opposite edges of the domain. 
To verify the new estimator $\hat{\eta}_{est}$ we employ
AVS-FE discretizations of the primal and dual problem with polynomial degrees of 
approximation of 2 and 3, respectively. 
The corresponding numerical results are illustrated in Table~\ref{tab:new_results_conv_diff}.
The effectivity indices 
show that the error estimator accurately measures the approximation error.  
\begin{table}[h]
\centering
\caption{\label{tab:new_results_conv_diff} Results for the convection-dominated diffusion problem \eqref{eq:conv_diff_ex} with QoI \eqref{eq:QoI_example_1} for uniform mesh refinements and $\Pe = 100$.}
\begin{tabular}{@{}lllll@{}}
\toprule
{Primal dofs \hspace{6mm}} & {$Q(u, \qq)-Q(u^h, \qq^h) $ \hspace{6mm}} & {Dual dofs \hspace{6mm}} & {$ \! \hat{\eta}_{est} \!$ \hspace{15mm}} & {$\hat{\eff} \!$} \\
\midrule \midrule

243 & 2.8825e-01 & 507 & 1.6140e-01 & 0.599   \\
867 & 1.7711e-01 & 1875 & 1.5010e-01 & 0.848   \\
3267 & 6.7393e-02 & 7203 & 6.5077e-02 & 0.966    \\
12675 & 1.3225e-02 & 28227 & 1.3158e-02 & 0.995   \\
49923 & 1.3918e-03 & 111747 & 1.3909e-03 & 0.999   \\
198147 & 1.0321e-04 & 444675 & 1.0321e-04 & 0.999   \\
 
\bottomrule
\end{tabular}
\end{table}
For initial, coarse, meshes the effectivity index may not be very close to unity, but is of the same order of magnitude as the
exact error. As the uniform meshes are further refined, the estimate converges 
to the exact error and delivers highly accurate predictions of the error with the effectivity index $\hat{\eff}$ very close to unity. 
\textcolor{black}{The rate of convergence of the error estimate in Table~\ref{tab:new_results_conv_diff} approaches the often observed superconvergence rate of 4 for bounded linear QoIs as reported by Giles and S{\"u}li in~\cite{giles2002adjoint}. }

The second QoI we consider is the average flux in the $x$-direction in the \textcolor{black}{same quadrant of the  
 unit square as before, i.e., $\omega = (0.5,1)\times(0.5,1)$. First, we consider the the derivative of the base variable $u$}:
\begin{equation}
\label{eq:QoI_example_3}
Q(u, \qq) = \frac{1}{|\omega|} \int_{\omega} \, \frac{\partial u}{\partial x}  \dx, 
\end{equation}
\textcolor{black}{and second, the flux variable:
\begin{equation}
\label{eq:QoI_example_3b}
Q(u, \qq) = \frac{1}{|\omega|} \int_{\omega} \, q_x  \dx. 
\end{equation}
} The polynomial degrees of approximation are now 1 and 2 for the primal and dual problem, respectively.
\begin{table}[h]
\centering
\caption{\label{tab:new_results_conv_diff_ex3} Results for the convection-dominated diffusion problem \eqref{eq:conv_diff_ex}, with QoI~\eqref{eq:QoI_example_3} for uniform mesh refinements and $\Pe = 100$.}
\begin{tabular}{@{}lllll@{}}
\toprule
{Primal dofs \hspace{6mm}} & {$Q(u, \qq)-Q(u^h, \qq^h) $ \hspace{6mm}} & {Dual dofs \hspace{6mm}} & {$ \! \hat{\eta}_{est} \!$ \hspace{15mm}} & {$\eff \!$} \\
\midrule \midrule


\textcolor{black}{867} & \textcolor{black}{-3.6294e-01} & \textcolor{black}{3267} & \textcolor{black}{1.9153e-01} & \textcolor{black}{-0.5277}   \\
\textcolor{black}{3267} & \textcolor{black}{-1.826e-01} & \textcolor{black}{12675} & \textcolor{black}{-2.1098e-03} & \textcolor{black}{0.0157}   \\
\textcolor{black}{12675} & \textcolor{black}{-6.8810e-02} & \textcolor{black}{49923} & \textcolor{black}{-4.3560e-02} & \textcolor{black}{0.6881}   \\
\textcolor{black}{49923} & \textcolor{black}{-1.7803e-02} & \textcolor{black}{198147} & \textcolor{black}{-1.1682e-02} & \textcolor{black}{0.9448}  \\

\bottomrule
\end{tabular}
\end{table}
\begin{table}[h]
\centering
\caption{\label{tab:new_results_conv_diff_ex3b} \textcolor{black}{Results for the convection-dominated diffusion problem \eqref{eq:conv_diff_ex}, with QoI~\eqref{eq:QoI_example_3b} for uniform mesh refinements and $\Pe = 100$.}}
\begin{tabular}{@{}lllll@{}}
\toprule
{\textcolor{black}{Primal dofs} \hspace{6mm}} & { \textcolor{black}{$Q(u, \qq)-Q(u^h, \qq^h) $} \hspace{6mm}} & { \textcolor{black}{Dual dofs} \hspace{6mm}} & { \textcolor{black}{$ \! \hat{\eta}_{est} \!$} \hspace{15mm}} & { \textcolor{black}{$\eff \!$} } \\
\midrule \midrule

\textcolor{black}{867} & \textcolor{black}{8.7745e-03} & \textcolor{black}{3267} & \textcolor{black}{1.1171e-02} & \textcolor{black}{1.3346}   \\
\textcolor{black}{3267} & \textcolor{black}{2.7558e-03} & \textcolor{black}{12675} & \textcolor{black}{3.8641e-03} & \textcolor{black}{1.4021}   \\
\textcolor{black}{12675} & \textcolor{black}{7.3292e-04} & \textcolor{black}{49923} & \textcolor{black}{8.4109e-04} & \textcolor{black}{1.1476}   \\
\textcolor{black}{49923} & \textcolor{black}{1.8478e-04} & \textcolor{black}{198147} & \textcolor{black}{1.8759e-04} & \textcolor{black}{1.0152}  \\

\bottomrule
\end{tabular}
\end{table}
As shown in Tables~\ref{tab:new_results_conv_diff_ex3} \textcolor{black}{and~\ref{tab:new_results_conv_diff_ex3b}},  the errors in a QoI 
in terms of a derivative are slightly higher than those of the preceding numerical 
verification \textcolor{black}{for the former case}, which is to be expected (see~\cite{CaloRomkesValseth2018,eirik2019thesis}). \textcolor{black}{The results in the latter are superior, again this is to be expected since the approximation of derivatives is generally less accurate.}
While for coarse meshes the estimate \textcolor{black}{in terms of the QoI with the derivative} does not accurately assess the
error, it does capture the right order of magnitude and improves significantly upon mesh refinements
as the effectivity indices $\hat{\eff}$ approaches unity. \textcolor{black}{However, the estimate in terms of the flux variable exhibits significant accuracy even for coarse meshes.}

\textcolor{black}{ The final numerical verification we consider for uniform mesh partitions is a QoI that
is the  average flux in the $x$-direction along the line segment on the left edge of the unit square, i.e., $\omega$ is the line segment from $y=0.5$ to $y=0.75$:
\begin{equation}
\label{eq:QoI_example_4b}
Q(u, \qq) = \frac{1}{|\omega|} \int_{\omega} \, q_x  \dx. 
\end{equation}
This type of QoI is particularly important in engineering design applications where fluxes or stresses and strains are critical design  parameters. Furthermore, for such quantities, classical methods often 
require enhancement techniques to achieve adequate accuracy such as the biharmonic smoothing introduced in~\cite{giles2002adjoint}. Furthermore, to show that the error estimate remains highly accurate when 
the error in the QoI becomes very small we pick $\Pe = 10$.
The polynomial degrees of approximation are 2 and 3 for the primal and dual problem, respectively.}

\begin{table}[h]
\centering
\caption{\label{tab:new_results_conv_diff_ex4b} \textcolor{black}{Results for the convection-dominated diffusion problem \eqref{eq:conv_diff_ex}, with QoI~\eqref{eq:QoI_example_4b} for uniform mesh refinements and $\Pe = 10$.}}
\begin{tabular}{@{}lllll@{}}
\toprule
{\textcolor{black}{Primal dofs} \hspace{6mm}} & { \textcolor{black}{$Q(u, \qq)-Q(u^h, \qq^h) $} \hspace{6mm}} & { \textcolor{black}{Dual dofs} \hspace{6mm}} & { \textcolor{black}{$ \! \hat{\eta}_{est} \!$} \hspace{15mm}} & { \textcolor{black}{$\eff \!$} } \\
\midrule \midrule

\textcolor{black}{867} & \textcolor{black}{-3.0174e-04} & \textcolor{black}{3267} & \textcolor{black}{-7.1583e-05} & \textcolor{black}{0.2372}   \\
\textcolor{black}{3267} & \textcolor{black}{-2.0855e-05} & \textcolor{black}{12675} & \textcolor{black}{-1.6377e-05} & \textcolor{black}{0.7853}   \\
\textcolor{black}{12675} & \textcolor{black}{-1.3401e-06} & \textcolor{black}{49923} & \textcolor{black}{-1.3090e-06} & \textcolor{black}{0.9768}   \\
\textcolor{black}{49923} & \textcolor{black}{-8.4335e-08} & \textcolor{black}{198147} & \textcolor{black}{-8.5814e-08} & \textcolor{black}{1.0175}  \\

\bottomrule
\end{tabular}
\end{table}
\textcolor{black}{The results in Table~\ref{tab:new_results_conv_diff_ex4b} show that the alternative 
approach is capable of estimating the error in terms of local QoIs pertaining to fluxes
 across boundaries. As the meshes become finer, the estimate becomes more accurate. However, for the coarsest mesh, i.e., the first row in Table~\ref{tab:new_results_conv_diff_ex4b}, the estimate is still 
 within an order of magnitude of the exact error. 
  As the error becomes very small, the new alternative method still provides accurate 
estimates without the need for additional enhancement techniques. }

\subsubsection{Non-Uniform Mesh}
\label{sec:non_uniform_mesh}

\begin{figure}[h]
\centering
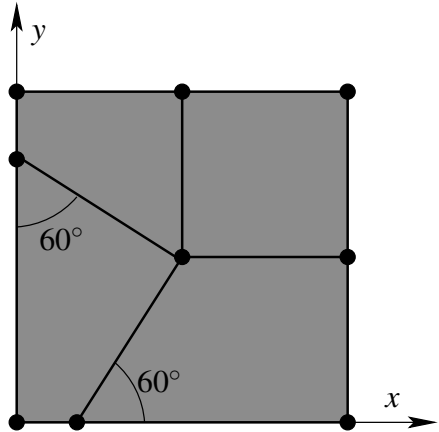 
\caption{Poorly constructed skewed mesh.}
\label{fig:skewed_mesh}
\end{figure}

So far, we have only considered rectangular uniform meshes. To provide a 
more realistic scenario, as encountered in engineering applications,
we consider a mesh in which the elements are skewed and the element edges do not align with the direction of 
the convection (see Figure~\ref{fig:skewed_mesh}).
We consider the same 
convection-dominated diffusion PDE~\eqref{eq:conv_diff_ex}, but now with
$\Pe=200$. 
The QoI is again the average of the solution $u$ in the region $\omega = (0.5,1)\times(0.5,1)$ (i.e., see~\eqref{eq:QoI_example_1}).

In Table~\ref{tab:skewed_results_conv_diff}, we list the results for the case
in which the primal degree of approximation is 2, and the dual degree of approximation
is 3. After each computation of the primal and dual solutions, all elements in the
mesh are uniformly refined. As in previous examples, the effectivity index is close to
unity, indicating that the estimator can be successfully employed for skewed meshes.
\begin{table}[h]
\centering
\caption{\label{tab:skewed_results_conv_diff} Results for the convection-dominated diffusion problem \eqref{eq:conv_diff_ex}, with QoI~\eqref{eq:QoI_example_1} for uniform mesh refinements of a skewed initial mesh partition and $\Pe = 200$.}
\begin{tabular}{@{}lllll@{}}
\toprule
{Primal dofs \hspace{6mm}} & {$Q(u, \qq)-Q(u^h, \qq^h) $ \hspace{6mm}} & {Dual dofs \hspace{6mm}} & {$ \! \hat{\eta}_{est} \!$ \hspace{15mm}} & {$\hat{\eff} \!$} \\
\midrule \midrule

867 & 1.9480e-01 & 1875 & 2.5397e-01 & 1.304   \\
3267 & 3.9104e-02 & 7203 & 4.6311e-02 & 1.184   \\
12675 & 4.0197e-03 & 28227 & 4.2989e-03 & 1.069   \\
49923 & 1.7641e-04 & 111747 & 1.7937e-04 & 1.017   \\
198147 & 5.0934e-06 & 444675 & 5.0988e-06 & 1.001  \\

\bottomrule
\end{tabular}
\end{table}

\subsubsection{Raviart-Thomas Approximation Of Fluxes}
\label{sec:Proper_Subspaces}
Until this point, we have used $C^0$ approximations for both trial 
variables, as our experience has shown this to yield good approximations~\cite{CaloRomkesValseth2018}.
Because $\qq \in \SHdivO$, the $C^0$ approximations have an overly restrictive
regularity.  
Commonly, in mixed FE methods, Raviart-Thomas rather than $\SCZO$ approximations are used. To show our approach also provides reliable estimates for such approximations 
 we now consider the case in which we use a Raviart-Thomas approximation~\cite{RaviartThomas,Girault1986,BrezziMixed} for the variables in $\SHdivO$ (i.e., $\qq$ and $\rr$).   
We again consider the case  where $\bb=\{1,1\}^T$, $\Pe=100$, and choose the same
QoI~\eqref{eq:QoI_example_1}, i.e., the average solution $u$ in
the top right quadrant of the unit square.
To approximate the error in the QoI, we now use tetrahedral elements in which $u^h,p^h$
are discretized with $C^0$ polynomials, while $\qq^h,\rr^h$ are
discretized using Raviart-Thomas bases. The initial mesh consists of two triangles which 
are refined uniformly after each computation using orders of approximation of 2 for the 
primal and 3 for the dual problems, respectively. 

For the same element partition, we achieve slightly higher accuracy at a slightly lower
number of degrees of freedom for  $C^0$ approximations versus Raviart-Thomas approximations for $\qq^h$, as evident from Tables~\ref{tab:new_results_conv_diff_ex5} and~\ref{tab:RT_results_conv_diff_ex5}. Comparison of the results in these tables also reveal that 
there is no significant difference between the two approximations in terms of the 
accuracy of the error estimates. 
Raviart-Thomas approximations are used in mixed FE methods as they result in stable FE approximations, \textcolor{black}{as well as consistency of the approximations}. Contrarily, $C^0$ approximations
for $\SHdivO$ variables cannot be employed in the same straightforward manner
for mixed FE methods and 
will lead to a violation of discrete \emph{inf-sup} conditions~\cite{BrezziMixed}. 
However, in the AVS-FE method, this \textcolor{black}{stability} problem is avoided by employing the DPG philosophy and
optimal test functions that ensure the discrete \emph{inf-sup} condition.
\begin{table}[h!]
\centering
\caption{\label{tab:new_results_conv_diff_ex5} Results for the convection-dominated diffusion problem \eqref{eq:conv_diff_ex}, with QoI~\eqref{eq:QoI_example_1} for uniform mesh refinements with a $C^0$ approximation for both variables and $\Pe = 100$. }
\scalebox{.8}{
\begin{tabular}{@{}llllll@{}}
\toprule
{Primal dofs \hspace{6mm}} & {$Q(u, \qq)-Q(u^h, \qq^h) $ \hspace{6mm}} & {$\! \norm{\qq-\qq^h}{\SLTO} \!$ \hspace{6mm}} & {Dual dofs \hspace{6mm}} & {$ \! \hat{\eta}_{est} \!$ \hspace{15mm}} & {$\hat{\eff} \!$} \\
\midrule \midrule

243 & 3.1381e-01 & 8.3021e-02 & 507 & 1.7312e-01  & 0.552    \\
867 & 1.9449e-01 & 5.7260e-02 & 1875 & 1.5716e-01 & 0.808   \\
3267 & 7.3123e-02 & 3.3955e-02 & 7203 & 7.2499e-02 & 0.991   \\
12675 & 1.3955e-02 & 1.4723e-02 & 28227 & 1.4085e-02 & 1.009    \\
49923 & 1.4397e-03 & 4.6769e-03 & 111747 & 1.4432e-03 & 1.002    \\

\bottomrule
\end{tabular} }
\end{table}
\begin{table}[h!]
\centering
\caption{\label{tab:RT_results_conv_diff_ex5} Results for the convection-dominated diffusion problem \eqref{eq:conv_diff_ex}, with QoI~\eqref{eq:QoI_example_1} for uniform mesh refinements with a Raviart-Thomas approximation of fluxes and $\Pe = 100$. }
\scalebox{.8}{
\begin{tabular}{@{}llllll@{}}
\toprule
{Primal dofs \hspace{6mm}} & {$Q(u, \qq)-Q(u^h, \qq^h) $ \hspace{6mm}} & {$\! \norm{\qq-\qq^h}{\SLTO} \!$ \hspace{6mm}} & {Dual dofs \hspace{6mm}} & {$ \! \hat{\eta}_{est} \!$ \hspace{15mm}} & {$\hat{\eff} \!$} \\
\midrule \midrule

257 & 3.5772e-01 & 9.0700e-02 & 529 & 1.8006e-01  & 0.503    \\
961 & 2.3563e-01 & 6.1567e-02 & 2017 & 1.7260e-01 & 0.733   \\
3713 & 9.9548e-02 & 3.7049e-02 & 7873 & 9.5166e-02 & 0.956   \\
14593 & 2.0993e-02 & 1.7127e-02 & 31105 & 2.1050e-02 & 1.003    \\
57857 & 2.2705e-03 & 5.7263e-03 & 123649 & 2.2741e-03 & 1.001    \\

\bottomrule
\end{tabular}}
\end{table}

\section{Adaptive Mesh Refinement}
\label{sec:adaptivity}
To demonstrate application of 
the new alternative error estimate~\eqref{eq:error_estimate} and the resulting 
error indicators~\eqref{eq:error_indicator} in an $h$-adaptive process, 
we use the same form of our model problem,
i.e., $\bb=\{1,1\}^T$ and $\Pe=100$. As the adaptive strategy for goal-oriented mesh refinement we use the method by Oden and Prudhomme~\cite{oden2001goal}, i.e.,
\begin{equation} \label{eq:ref_crit_goal}
\text{if } \,\,\, \frac{\ds | \epsilon_m | }{\ds \text{max}_{\Kep} | \epsilon_m |} > \delta, \,\,\,
\text{then refine element} \,\,m,
\end{equation}
where $\delta$ is the tolerance for refinement, i.e., $0 < \delta < 1$. In the following numerical verification, we pick $\delta = 0.5$.
The QoI we consider is again the average of $u$ in the upper right quadrant (see~\eqref{eq:QoI_example_1}).
The primal problem is approximated using $C^0$ continuous polynomials of degree 2, whereas the dual problem is approximated using $p+1=3$.

\begin{figure}[h]
{\centering
\hspace{1.1in}  \includegraphics[width=0.5\textwidth]{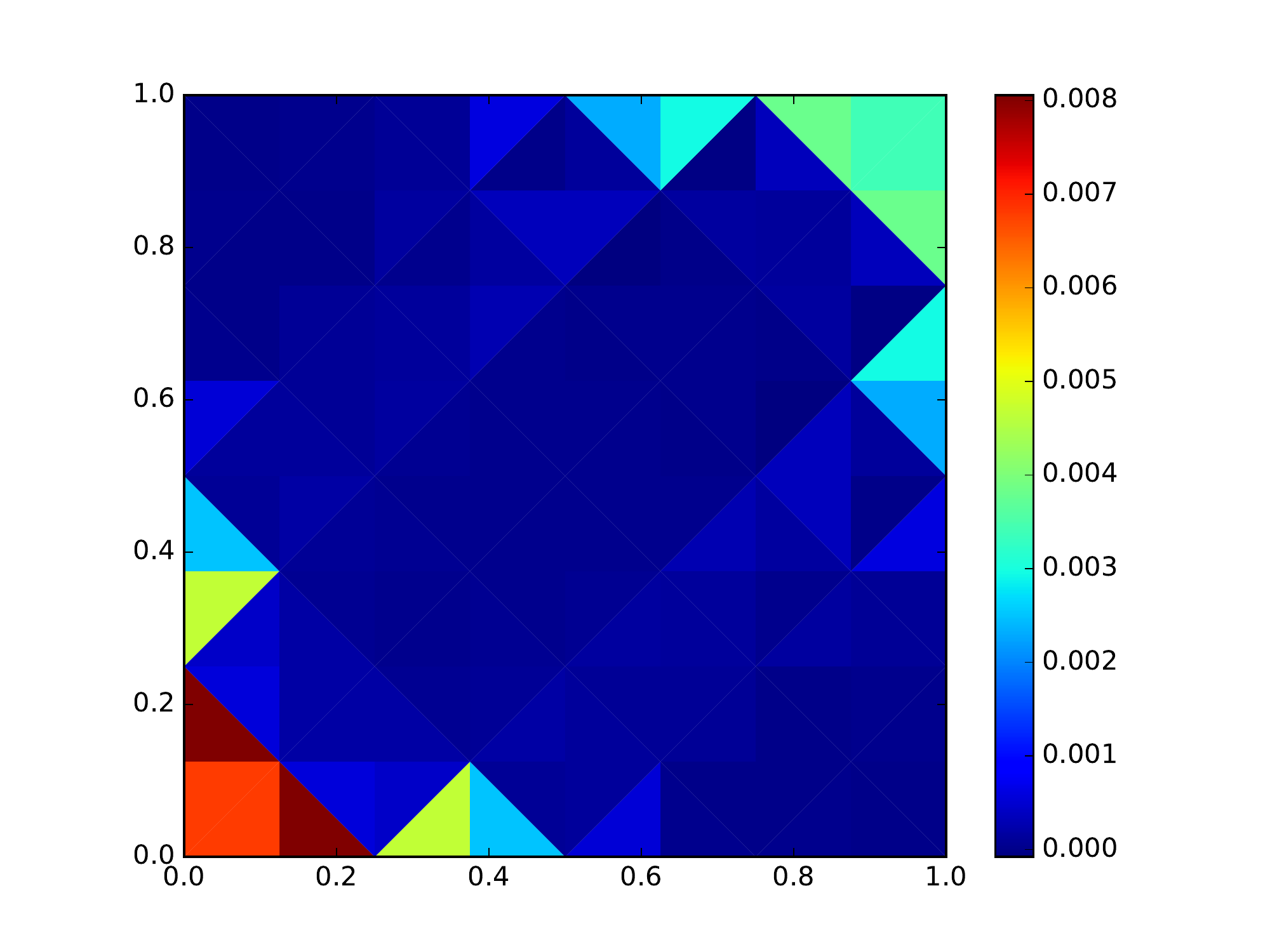}}
 \caption{\label{fig:local_err_ind_s2_goal} Element-wise distribution of the error indicators $\epsilon_m $ \eqref{eq:error_indicator} on a coarse mesh. }
\end{figure}

The initial mesh consists of 2 triangular elements and is too coarse 
to resolve the boundary layers in both primal and dual solutions leading to poor error indicators. This effect is shown in Figure \ref{fig:local_err_ind_s2_goal}, where the error indicators are largest in the 
corner of the dual boundary layer which would result in mesh refinements at the 'wrong' location.  To avoid initial mesh 
refinements that are poorly suited to reduce the error in the QoI, we
initially perform uniform mesh refinements until the error estimate $\hat{\eta}_{est}$ begins to decrease and indicate that the error indicators $\epsilon_m$~\eqref{eq:error_indicator} have become reliable.

\begin{figure}[h]
\subfigure[ \label{fig:2d_goal_sol} Solution $u^h$. ]{\centering
 \includegraphics[width=0.45\textwidth]{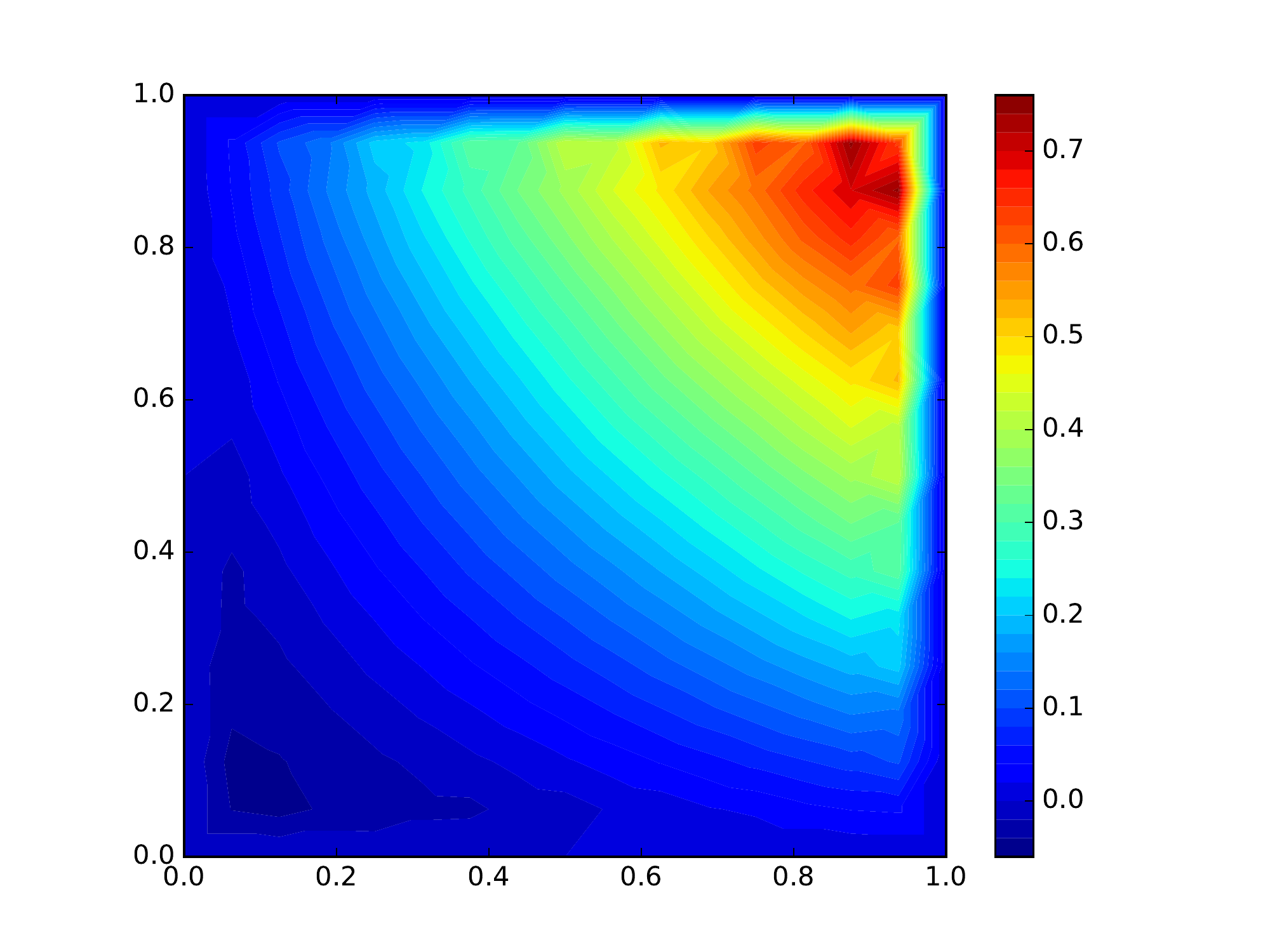}}
\hspace{0.1in}    \subfigure[ \label{fig:2d_goal_sol_indicators} Error indicator $\epsilon_m$.]{\centering
 \includegraphics[width=0.45\textwidth]{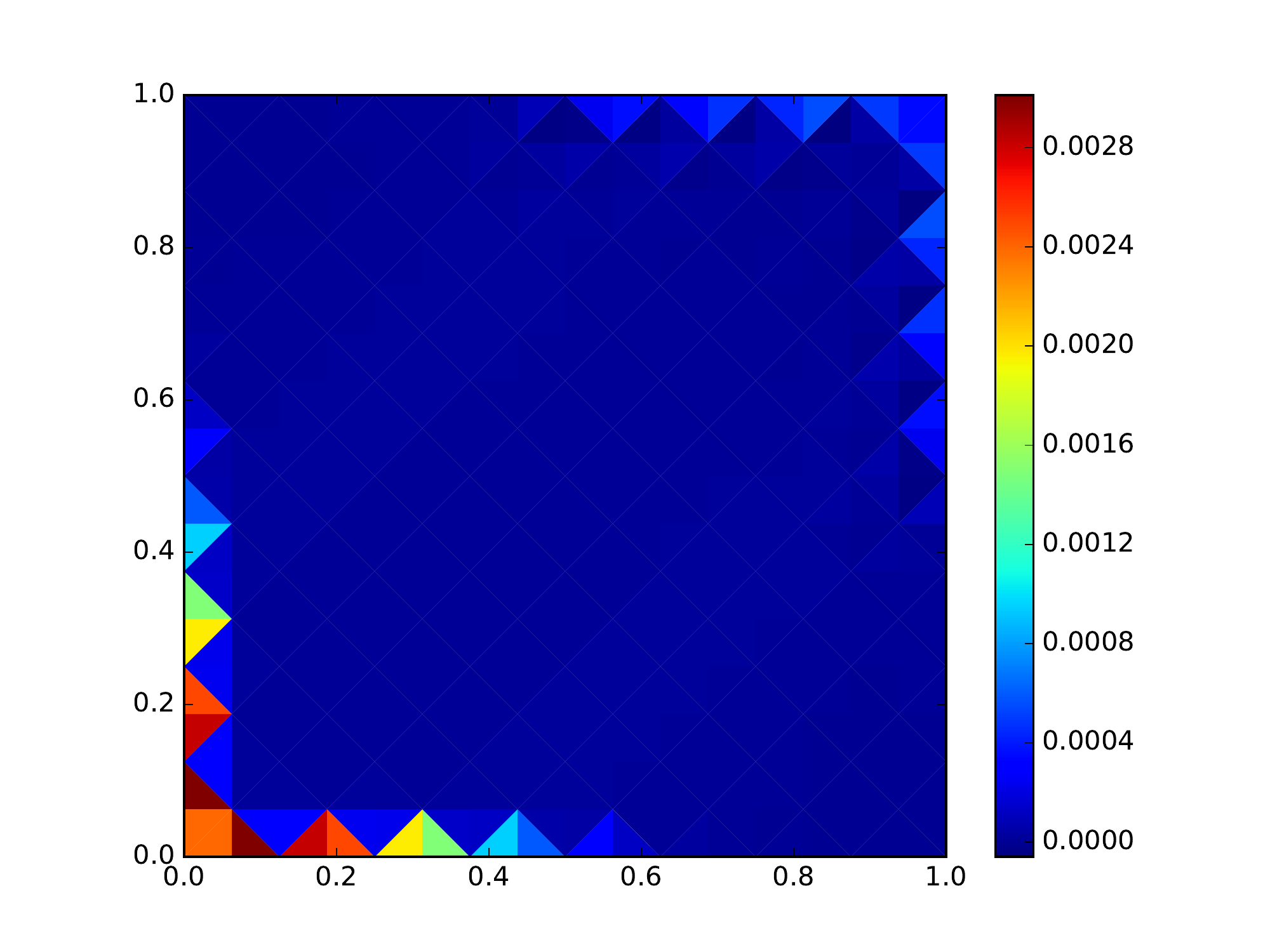}}
  \caption{\label{fig:initial_goal_step} Solution on a uniform $16\times16$ mesh.}
\end{figure}

\begin{figure}[h]
\subfigure[ \label{fig:2d_goal_sol_indicators_9} Step 9. ]{\centering
 \includegraphics[width=0.45\textwidth]{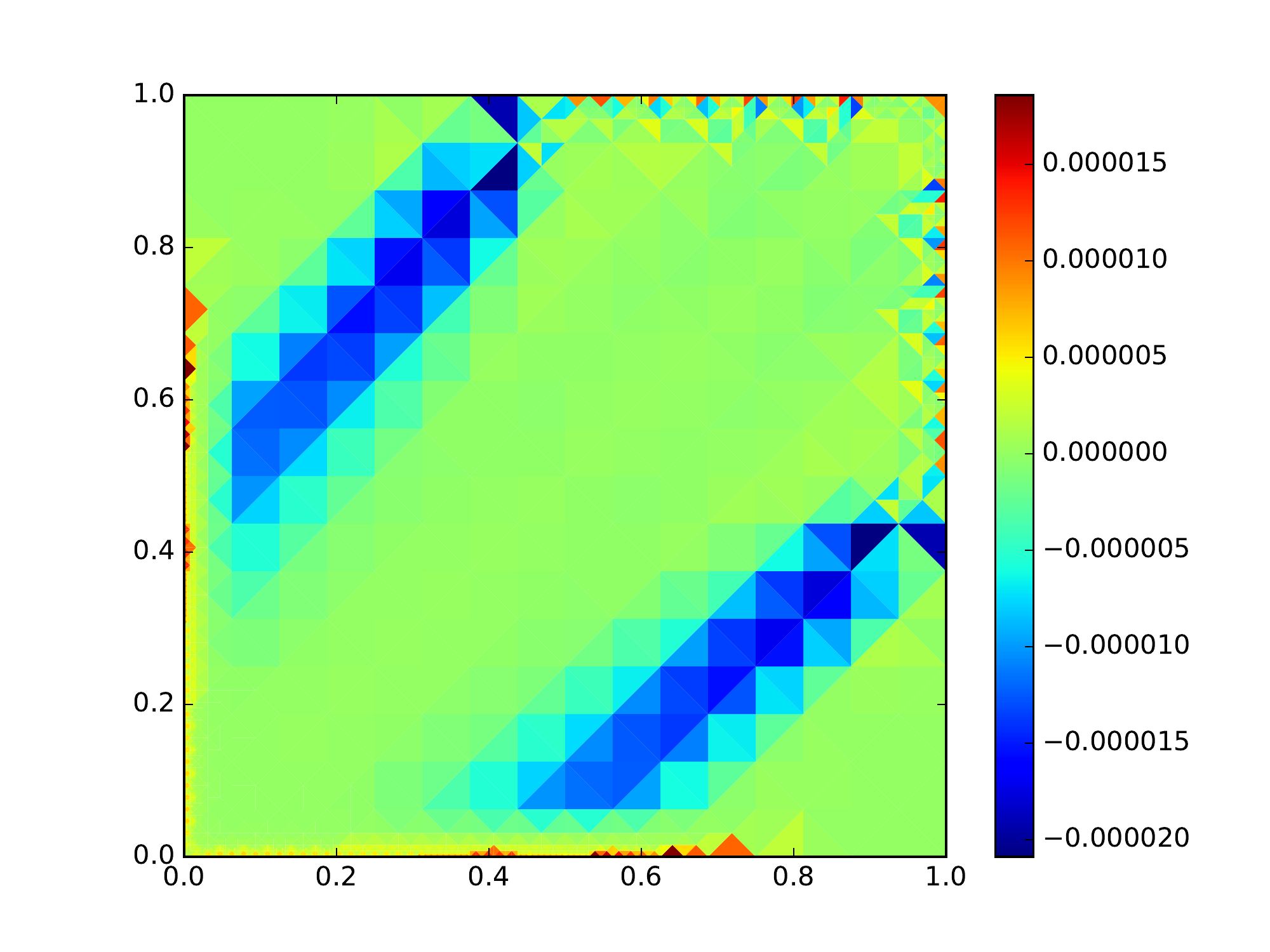}}
\hspace{0.1in}    \subfigure[ \label{fig:2d_goal_sol_indicators_16} Step 18.]{\centering
 \includegraphics[width=0.45\textwidth]{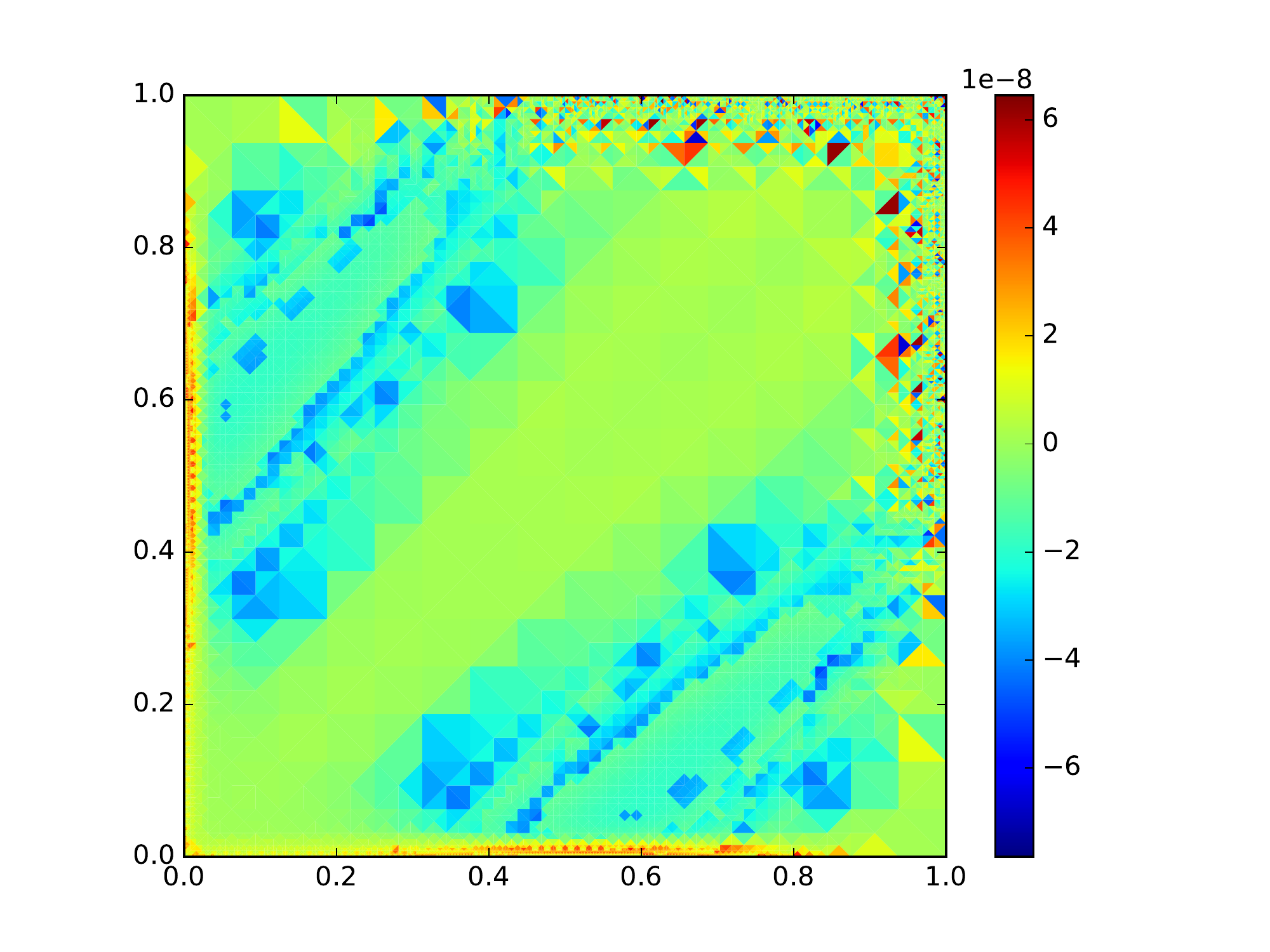}}
  \caption{\label{fig:goal_step9_16} Error indicators $\epsilon_m$.}
\end{figure}

\begin{figure}[h]
{\centering
\hspace{1.1in}  \includegraphics[width=0.5\textwidth]{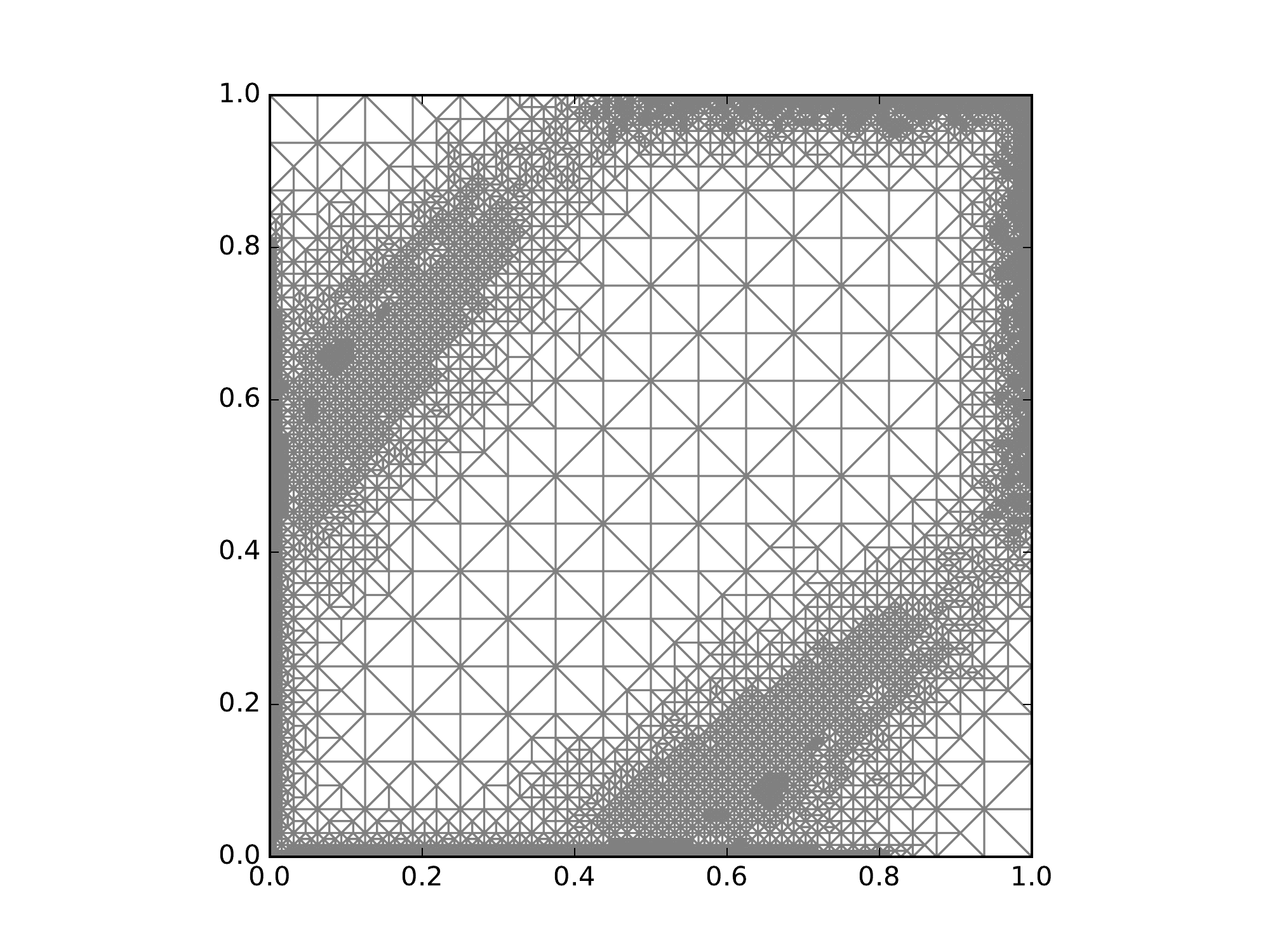}}
 \caption{\label{fig:goal_step_16_mesh} Final mesh of the goal-oriented $h-$adaptive refinements. }
\end{figure}

\begin{figure}[h]
\subfigure[ \label{fig:2d_goal_sol_converg} $\hat{\eta}_{est}$ and $\norm{u-u^h}{\SLTO}$.]{\centering
 \includegraphics[width=0.5\textwidth]{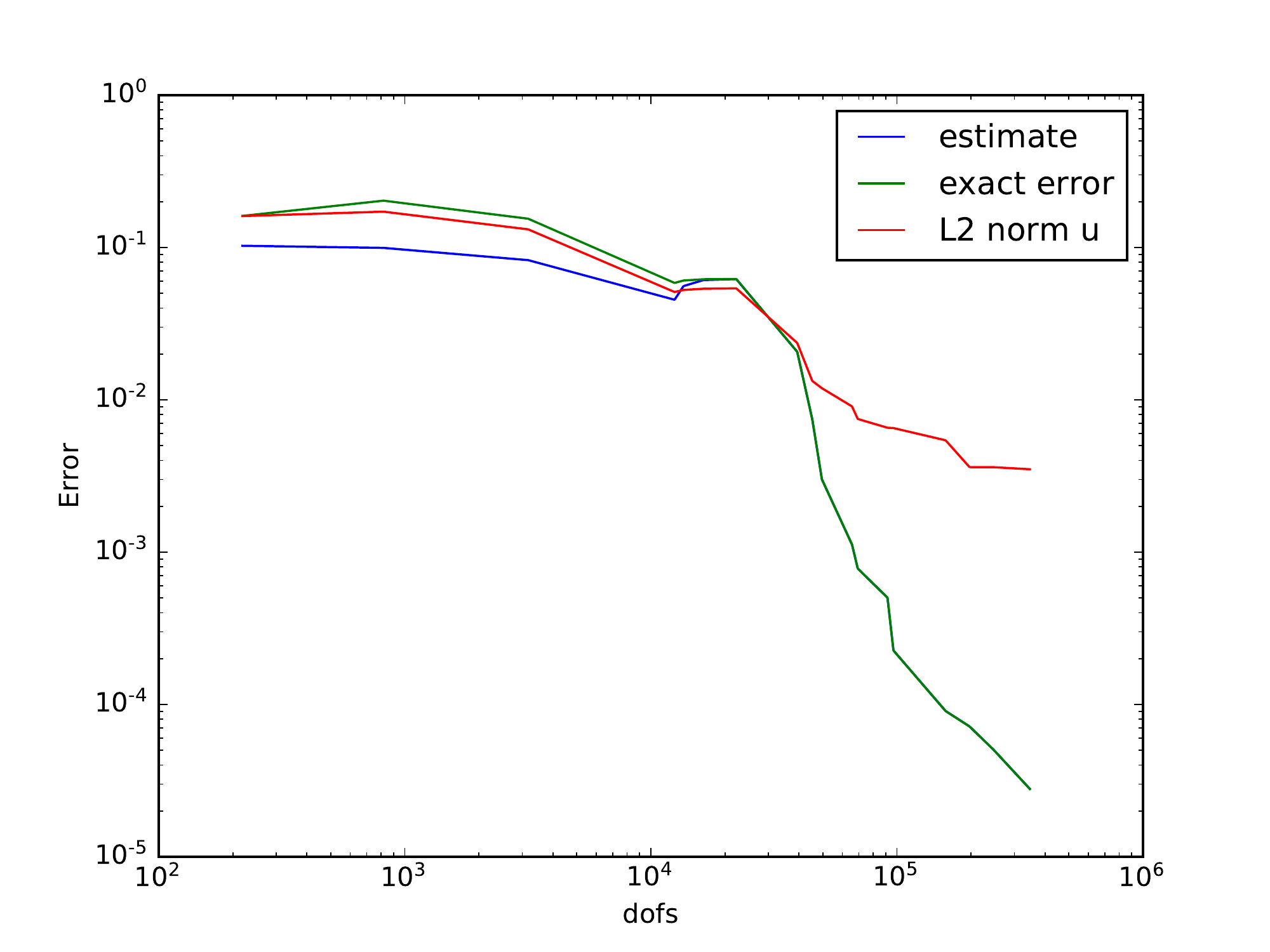}} \hspace{0.1in} 
\subfigure[ \label{fig:2d_goal_sol_eff_ind} Effectivity index $\hat{\eff}$. ]{\centering
 \includegraphics[width=0.5\textwidth]{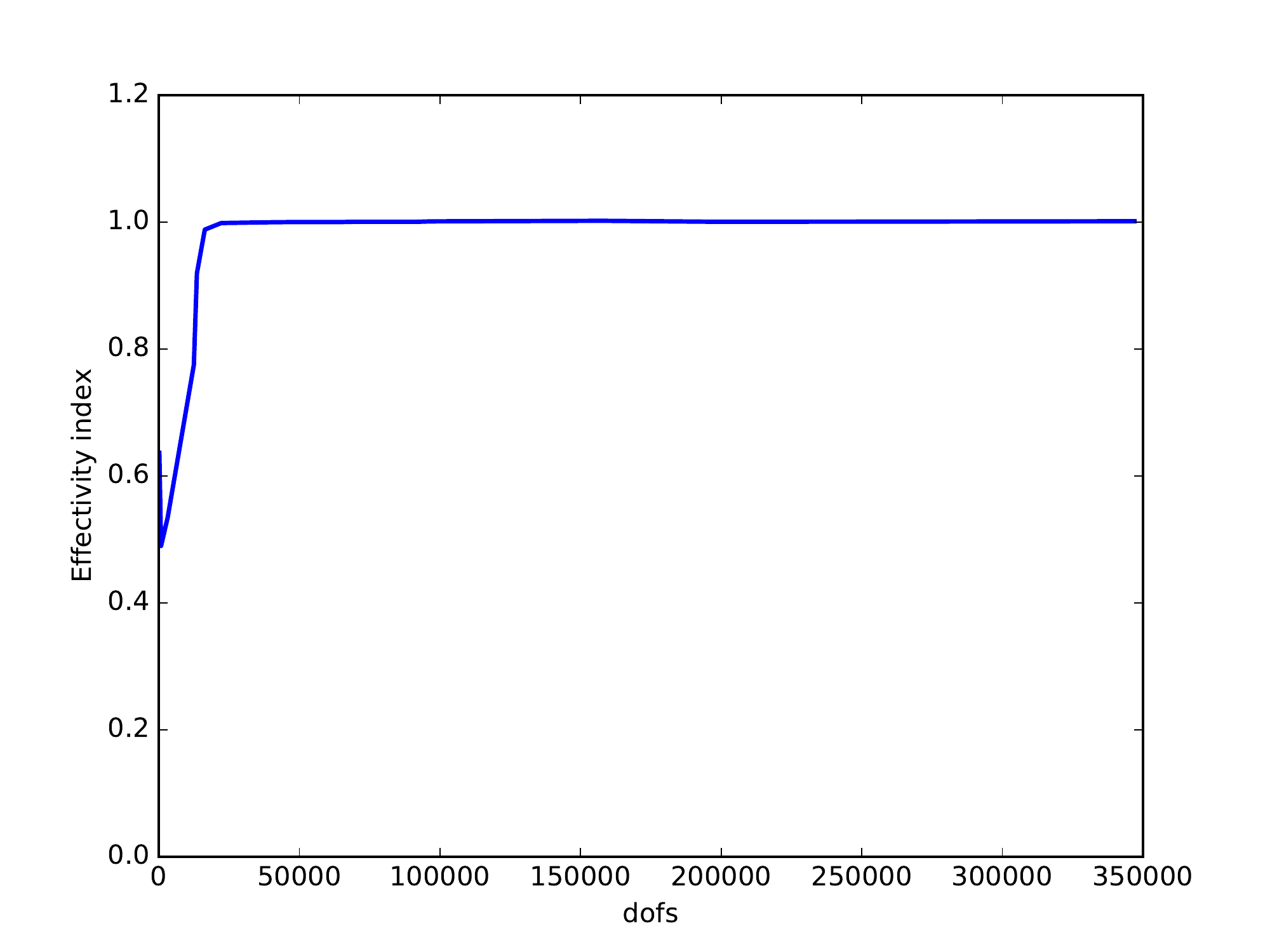}}
   
  \caption{\label{fig:2d_goal_sol_converg_results} Convergence history for goal-oriented adaptive mesh refinement.}
\end{figure}

This can for example be seen in Figure \ref{fig:initial_goal_step}, where we show the  solution and element-wise
error indicators on a uniform mesh of $16\times16$ elements, the last mesh that has been uniformly refined and $\hat{\eta}_{est}$ started to decrease.  Here, we see that  indicators in corner of the primal boundary layer are now of a magnitude that result in local mesh refinements in the right locations. 
In Figure \ref{fig:goal_step9_16}, we show the error indicators for an intermediate (step 9) and the final (step 18) step of the adaptive process. 
In both cases, the mesh has been refined such that the boundary layers in both primal and dual problems are sufficiently resolved to yield error indicators that are
highest in the region of the QoI.  The 
corresponding final adapted mesh is shown in Figure \ref{fig:goal_step_16_mesh}. 
As expected from the current choice of QoI, the mesh refinements have been focused near 
the primal boundary layer and the QoI.
Lastly, we present the convergence history of the error estimator $\hat{\eta}_{est}$ 
and the effectivity index $\hat{\eff}$ in Figure~\ref{fig:2d_goal_sol_converg_results}. The plot of estimated error and  $\norm{u-u^h}{\SLTO}$ in Figure
\ref{fig:2d_goal_sol_converg} shows that while the adaptive process ensures a small
error in the QoI, the global error in the $\SLTO$ norm is several orders of magnitude
larger as expected since the adaptive process is targeting to reduce the error in the QoI rather than $\norm{u-u^h}{\SLTO}$.
The effectivity index shown in Figure~\ref{fig:2d_goal_sol_eff_ind} demonstrates that the 
proposed alternative approach to goal-oriented error estimation delivers highly 
accurate estimates even when the error becomes very small.

\begin{figure}[h]
\subfigure[ \label{fig:2d_goal_sol_converg2}\textcolor{black}{ $\hat{\eta}_{est}$ and $\norm{u-u^h}{\SLTO}$.}]{\centering
 \includegraphics[width=0.5\textwidth]{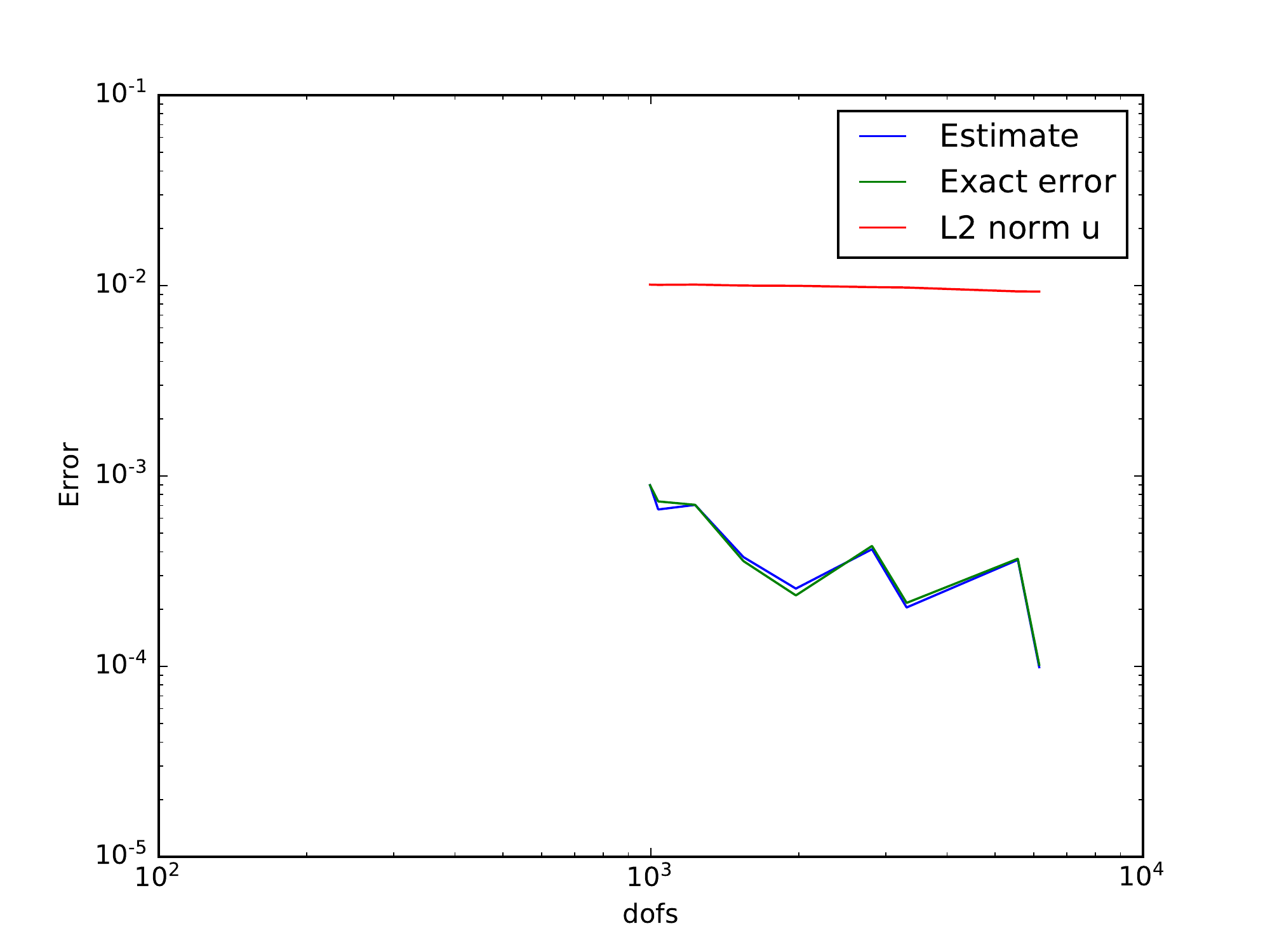}} \hspace{0.1in} 
\subfigure[ \label{fig:2d_goal_sol_eff_ind2} \textcolor{black}{ Effectivity index $\hat{\eff}$. }]{\centering
 \includegraphics[width=0.5\textwidth]{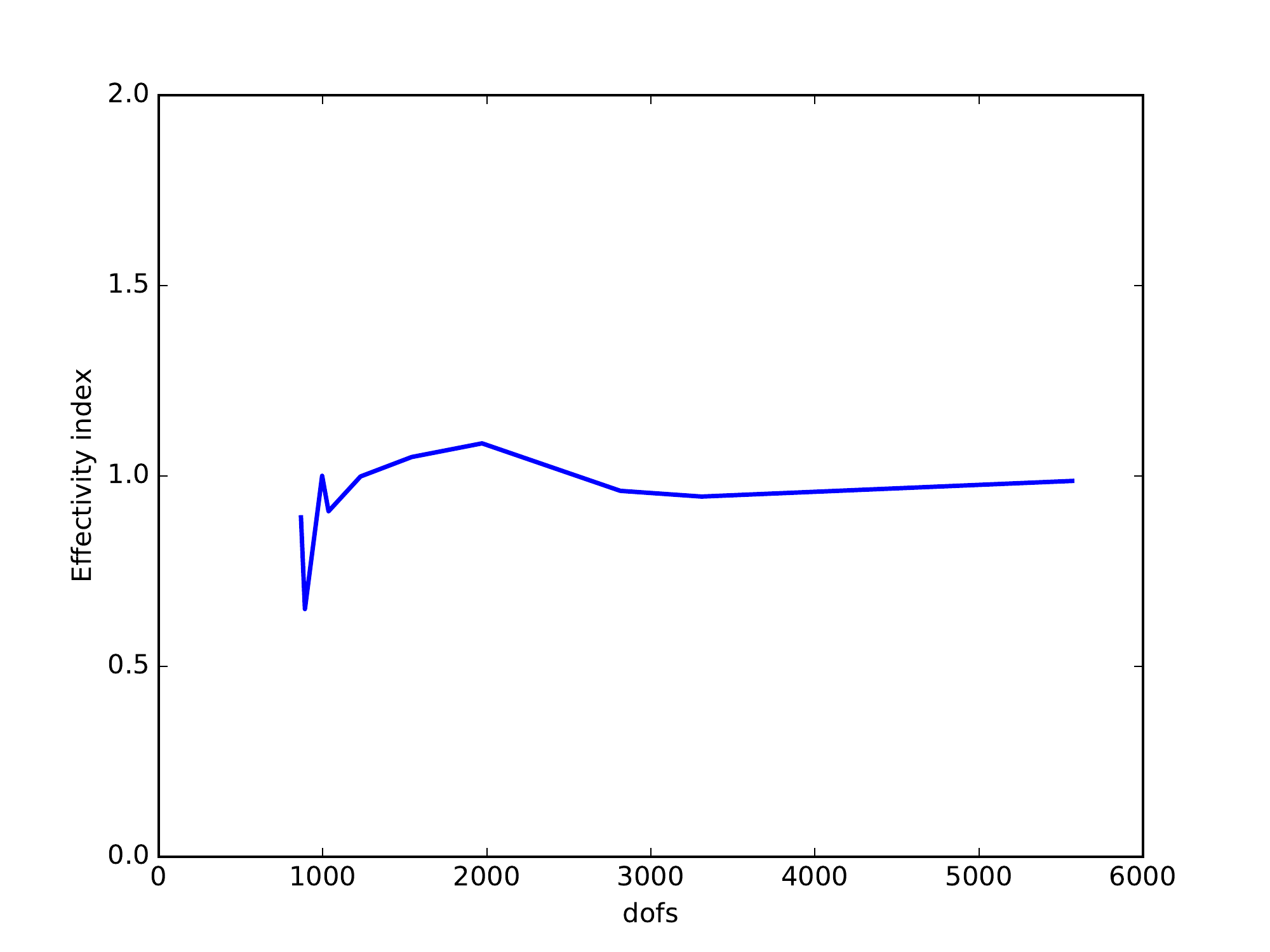}}
  \caption{\label{fig:2d_goal_sol_converg_results2} \textcolor{black}{Convergence history for goal-oriented adaptive mesh refinement a flux QoI~\eqref{eq:QoI_example_4b}.}}
\end{figure}

\begin{figure}[h]
{\centering
\hspace{1.1in}  \includegraphics[width=0.5\textwidth]{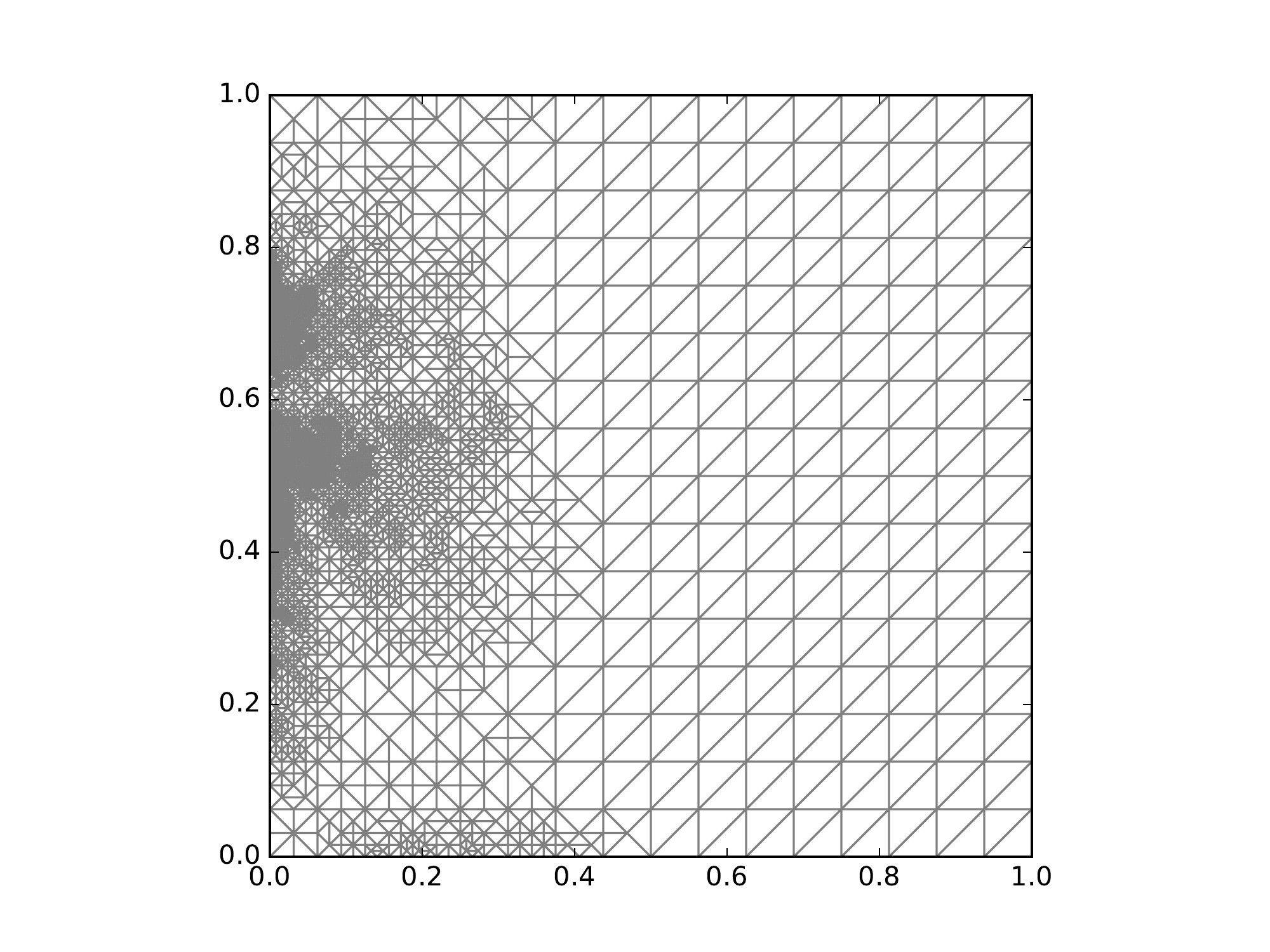}}
 \caption{\label{fig:goal_step_10_mesh} \textcolor{black}{Final mesh of the goal-oriented $h-$adaptive refinements using a flux QoI~\eqref{eq:QoI_example_4b}. }}
\end{figure}

\textcolor{black}{As a final numerical verification we consider the flux $q_x$ across the line between $y=0.5$ and $y=0.75$ on the left edge of the unit square in~\eqref{eq:QoI_example_4b}, and we keep the same physical parameters as the corresponding verification, i.e., $\bb=\{1,1\}^T$ and $\Pe=10$.  
The primal problem is approximated using $C^0$ continuous polynomials of degree 2, whereas the dual problem is approximated using the same type polynomials at degree $p+1=3$. The initial mesh consists of 
512 triangle elements,  we employ the adaptive strategy as in~\eqref{eq:ref_crit_goal} with identical 
tolerance as well, and we perform 10 mesh refinements here. 
This initial mesh is chosen to be sufficiently fine to adequately resolve the boundary layers to 
reduce pollution effects.
In Figure~\ref{fig:2d_goal_sol_converg2}, the 
convergence history of the error indicator as well as the global $L^2$ error in $u$. The corresponding effectivity index  in Figure~\ref{fig:2d_goal_sol_eff_ind2} shows that the estimate remains highly accurate during the adaptive process. Lastly, the final adapted 
mesh shown in Figure~\ref{fig:goal_step_10_mesh} show that the refinements are focused near the 
QoI as expected for the current refinement criterion.}

\section{Conclusions}
\label{sec:conclusions}
We have presented goal-oriented \emph{a posteriori} error estimates for
the AVS-FE method. This method is a hybrid Petrov-Galerkin method which uses classical $\SCZO$ or Raviart-Thomas
 FE trial basis functions, while the test space consists of functions that are 
discontinuous across element edges. The broken topology of the test space allows us 
to employ the DPG philosophy and compute optimal test functions element-by-element, i.e., completely locally.
In an effort to derive \emph{a posteriori} error estimates of the AVS-FE computations 
we have introduced two types goal-oriented error estimates. The first estimate follows 
the classical approach of~\cite{becker_rannacher_2001} where, by duality,
the dual solution is sought in the test space $V$, which in the case of the AVS-FE method is 
a broken Hilbert space. However, we show that through numerical verifications of the 
classical Laplace BVP that 
this approach yields error estimates with poor accuracy.
To resolve this, we introduce a second estimate based on consideration of the PDE
that governs the dual solution. The estimate is then established by computing 
$\SCZO$ or Raviart-Thomas AVS-FE approximations of this PDE. Numerical verifications involving pure diffusion as well as convection-dominated diffusion
problems show that the new alternative error estimate is capable of accurately 
predicting errors for different QoIs and mesh partitions. 

In order to employ the new \emph{a posteriori}  error estimation methodology in mesh 
adaptivity we also here derived error indicators to guide any $h-$adaptive process.
The error indicators essentially are the element-wise restriction of the residual 
operator~\eqref{eq:error_indicator}.
Numerical verifications show that when the error indicators used with classical 
refinement strategies~\cite{oden2001goal}, they lead to a mesh adaptive process
able to reduce the error in the QoI within a defined accuracy while at the same 
time delivering highly accurate predictions of the error in the QoI even when the 
error is small.
In a previous paper~\cite{CaloRomkesValseth2018}, the results presented are all computed 
using $C^0$ approximations for both trial variables. 
Here, we also presented 
 results in which the fluxes are computed by using
Raviart-Thomas approximations. These results indicate that  $C^0$ approximations
yield results that are of the same quality in terms of both error estimation, and 
result in slightly higher accuracy for the flux variable at a slightly lower number of degrees of freedom. Hence, the use of $C^0$ 
approximations for both variables remains  attractive due to its lower
computational cost, and ease of implementation in existing FE software.

\textcolor{black}{Note that the domains considered in the verifications in Section~\ref{sec:Proper_Subspaces} are convex, it is likely that for non-convex domains the consistency of Raviart-Thomas approximations will be preferable over our $C^0$ approximations. This investigation is postponed to future research efforts.
The poor performance of the classical method reported in Section~\ref{sec:class_estimates} appears to be related to regularity of the dual solution as suggested by DPG$^*$ literature. Hence, in future efforts we will pursue analyses of the dual problem within a framework similar to the DPG$^*$ method 
to fully understand the intricacies of the dual solution.}

In a forthcoming paper, we intend to extend the AVS-FE method and the new alternative error estimates 
to other problems such as the nonlinear Cahn-Hilliard equation, \textcolor{black}{as well as alternative error indicators as proposed in, e.g.,~\cite{darrigrand2018goal}}.

\section*{Acknowledgements}
The numerical verifications presented in Section~\ref{sec:Proper_Subspaces} were 
computed using the computational framework Firedrake~\cite{rathgeber2017firedrake} and the $h-$adaptive 
process presented in Section~\ref{sec:adaptivity} was implemented in the computational framework FEniCS~\cite{alnaes2015fenics}.
\newline
This work has been supported by the \textcolor{black}{United States National Science Foundation -} NSF CBET Program,
under  NSF Grant  titled \emph{Sustainable System for Mineral Beneficiation}, NSF Grant
No. 1805550.

\bibliographystyle{elsarticle-num}
 \bibliography{references_eirik}
\end{document}